\documentclass[11pt,a4paper]{article}
\date{}
\makeatletter
\@addtoreset{equation}{section}
\makeatother

\usepackage[colorlinks,
            linkcolor=blue,
            anchorcolor=blue,
            citecolor=blue]{hyperref}
\usepackage{hyperref}
\hypersetup{hypertex=true,
            colorlinks=true,
            linkcolor=blue,
            anchorcolor=blue,
            citecolor=blue}
\usepackage{hyperref,amsmath,amssymb,amscd}
\usepackage{amssymb,amsmath,enumerate}
\usepackage{indentfirst}
\usepackage{url}
\usepackage[capitalise]{cleveref}
\usepackage{amsthm}
\usepackage{amsmath}
\newtheorem{thm}{Theorem}[section]
\newtheorem{cor}[thm]{Corollary}
\newtheorem{lem}[thm]{Lemma}

\theoremstyle{definition}

\newtheorem{rem}[thm]{\bf Remark}
\allowdisplaybreaks[4]

\usepackage{enumitem}
\usepackage{extarrows}
\usepackage{amsthm}
\usepackage{setspace}
\hypersetup{urlcolor=blue, citecolor=red}
\usepackage{graphicx}
\usepackage{tikz}

\begin{document}

\title{\bf
The second order Caffarelli-Kohn-Nirenberg identities and inequalities\footnote{
Supported by National Natural Science Foundation of China (No. 12371120).}}
\author{{Xiao-Ping Chen,\ \ Chun-Lei Tang\footnote{Corresponding author.\newline
\indent\,\,\, \emph{E-mail address:} tangcl@swu.edu.cn (C.-L. Tang); xpchen\_maths@163.com (X.-P. Chen).}}\\
{\small \emph{School of Mathematics and Statistics, Southwest University,  Chongqing {\rm400715},}}\\
{\small \emph{People's Republic of China}}}
\maketitle
\baselineskip 17pt

\noindent {\bf Abstract}:\ This paper focuses on optimal constants and optimizers of the second order Caffarelli-Kohn-Nirenberg inequalities. Firstly, we aim to study optimal constants and optimizers for the following second order Caffarelli-Kohn-Nirenberg inequality in radial space: let $N\ge1$, $t\ge p>1$,
\begin{equation}\label{0.1}
\left(\int_{\mathbb{R}^N}
\frac{|\Delta u|^p}{|x|^{p\alpha}}
\mathrm{d}x\right)^{\frac{1}{p}}
\left[\int_{\mathbb{R}^N}
\frac{\left|\nabla u\right|^{\frac{p(t-1)}{p-1}}}
{|x|^{\frac{p(t-1)}{p-1}\beta}}
\mathrm{d}x\right]^{\frac{p-1}{p}}
\ge C(N,p,t,\alpha,\beta)
\int_{\mathbb{R}^N}
\frac{\left|\nabla u\right|^t}{|x|^{t\gamma}}
\mathrm{d}x.
\end{equation}
Secondly, we establish second order  $L^p$-Caffarelli-Kohn-Nirenberg identities, and obtain optimal constants and optimizers of the second order $L^p$-Caffarelli-Kohn-Nirenberg inequalities (i.e., $p=t$ in \eqref{0.1}) in general space. Lastly, under some more general assumptions, we consider the optimal weighted second order Heisenberg Uncertainty Principles, which complements the recent work [``The sharp second order Caffareli-Kohn-Nirenberg inequality and stability estimates for the sharp second order uncertainty principle'', 2022, arXiv:2102.01425].

This paper's main novelty lies in the fact that we research the optimal versions of the second order Caffarelli-Kohn-Nirenberg inequalities \eqref{0.1} in radial space or in general space, and also establish the second order $L^p$-Caffarelli-Kohn-Nirenberg identities.

\vspace{0.25em}


\vspace{0.25em}


\section{Introduction}

\noindent In this paper, we aim to study optimal constants and optimizers of the second order  Caffarelli-Kohn-Nirenberg (CKN for short) inequalities and the weighted second order Heisenberg Uncertainty Principles (HUP for short). Our goals of this paper are threefold.
\begin{enumerate}
[itemsep=0pt, topsep=2pt, parsep=0pt]

\item[($a$)] Firstly, optimal constants and optimizers will be considered for a more general subfamily of the second order CKN inequalities of the following form: let $N\ge1$ and $t\ge p>1$,
\begin{equation}\label{1.1}
\left(\int_{\mathbb{R}^N}
\frac{|\Delta u|^p}{|x|^{p\alpha}}
\mathrm{d}x\right)^{\frac{1}{p}}
\left[\int_{\mathbb{R}^N}
\frac{\left|\nabla u\right|^{\frac{p(t-1)}{p-1}}}
{|x|^{\frac{p(t-1)}{p-1}\beta}}
\mathrm{d}x\right]^{\frac{p-1}{p}}
\ge C_1(N,p,t,\alpha,\beta)
\int_{\mathbb{R}^N}
\frac{\left|\nabla u\right|^t}{|x|^{t\gamma}}
\mathrm{d}x,
\end{equation}
for all radial functions $u\in\mathcal{C}_0^\infty(\mathbb{R}^N
\setminus\{\mathbf{0}\})$, and the constant $C_1(N,p,t,\alpha,\beta)$ of \eqref{1.1} is optimal (see Theorem \ref{thm-2.1}). The main novelty in this case is that we consider the optimal second order CKN inequalities when $p\neq t$ and $p\neq2$.

\item[($b$)] Secondly, we establish second order $L^p$-CKN identities (see Theorems \ref{thm-2.3} and \ref{thm-2.4}), based on them, we also analyze the optimal second order $L^p$-CKN inequality
\begin{equation}\label{1.2}
\left(\int_{\mathbb{R}^N}\frac{|\Delta u|^p}
{|x|^{p\alpha}}\mathrm{d}x\right)^{\frac{1}{p}}
\left(\int_{\mathbb{R}^N}\frac{|\nabla u|^p}
{|x|^{p\beta}}\mathrm{d}x\right)^{\frac{p-1}{p}}
\ge C_2(N,p,\alpha,\beta)
\int_{\mathbb{R}^N}
\frac{|\nabla u|^p}{|x|^{\alpha+(p-1)\beta+1}}\mathrm{d}x.
\end{equation}
for all $u\in\mathcal{C}_0^\infty(\mathbb{R}^N
\setminus\{\mathbf{0}\})$, and the constant  $C_2(N,p,\alpha,\beta)$ of \eqref{1.2} is optimal (see Theorems \ref{thm-2.5} and \ref{thm-2.6}). Here we would like to point out that, according to the second order $L^p$-CKN identities and a function $\mathcal{K}_p$ (defined in \eqref{2.3} below), radial assumption about $u$ is not required in this case.

\item[($c$)] Lastly, under more general conditions, we study the optimal weighted second order HUP
\begin{equation}\label{1.3}
\left(\int_{\mathbb{R}^N}
\frac{|\Delta u|^2}{|x|^{2\alpha}}\mathrm{d}x
\right)^{\frac{1}{2}}
\left(\int_{\mathbb{R}^N}
|x|^{2\alpha+2}|\nabla u|^2\mathrm{d}x
\right)^{\frac{1}{2}}
\ge C_3(N,\alpha)
\int_{\mathbb{R}^N}|\nabla u|^2\mathrm{d}x,
\end{equation}
for all $u\in\mathcal{C}_0^\infty(\mathbb{R}^N
\setminus\{\mathbf{0}\})$, and the constant  $C_3(N,\alpha)$ of \eqref{1.3} is optimal (see Theorem \ref{thm-2.8}). Observe that this result can be considered as a supplement to that of \cite{Duong21}.

\end{enumerate}

We review the research status of CKN inequalities and some related inequalities in Sect. \ref{1.1}. Some related results about optimal constants and optimizers of the first order and second order CKN inequalities will be presented in Sect. \ref{1.2} and Sect. \ref{1.3}, respectively.

\subsection{Overview of Caffarelli-Kohn-Nirenberg inequalities}\label{1.1}

\noindent In 1984, Caffarelli, Kohn and Nirenberg  \cite{Caffarelli84} firstly introduced Caffarelli-Kohn-Nirenberg inequalities: let $N\ge1$, parameters $p,q,t,\alpha,\beta,\sigma$ and $\delta$ satisfy $p,q\ge1$, $t>0$, $0\le\delta\le1$,
$$
N-p\alpha>0,
\
N-q\beta>0,
\
N-t\gamma>0
$$
and
$$
\gamma=\delta\sigma+(1-\delta)\beta.
$$
Then for all $u\in\mathcal{C}_0^\infty(\mathbb{R}^N
\setminus\{\mathbf{0}\})$, there is a constant $C>0$ satisfying
\begin{equation}\label{1.4}
\left(\int_{\mathbb{R}^N}
\frac{|\nabla u|^p}
{|x|^{p\alpha}}\mathrm{d}x
\right)^{\frac{\delta}{p}}
\left(\int_{\mathbb{R}^N}
\frac{|u|^q}{|x|^{q\beta}}
\mathrm{d}x\right)^{\frac{1-\delta}{q}}
\geq C
\left(\int_{\mathbb{R}^N}
\frac{|u|^t}{|x|^{t\gamma}}
\mathrm{d}x\right)^{\frac{1}{t}}
\end{equation}
if and only if
$$
\frac{1}{t}-\frac{\gamma}{N}
=\delta\left(\frac{1}{p}-\frac{\alpha+1}{N}\right)
+(1-\delta)\left(\frac{1}{q}-\frac{\beta}{N}\right),
$$
(dimensional balance),
$$
\alpha-\sigma\ge0,
\
\mathrm{if}
\
\delta>0,
$$
and
$$
\alpha-\sigma\le1,
\
\mathrm{if}
\
\delta>0
\
\mathrm{and}
\
\frac{1}{p}-\frac{\alpha+1}{N}
=\frac{1}{t}-\frac{\gamma}{N}.
$$
It is well-known that there are many important inequalities included in the CKN inequality \eqref{1.4}, such as Sobolev inequalities ($\delta=1$, $\alpha=\gamma=0$), Gagliardo-Nirenberg inequalities ($0<\delta<1$, $\alpha=\beta=\gamma=0$), Hardy inequalities ($\delta=1$, $\alpha=0$, $p=t$), Hardy-Sobolev inequalities ($\delta=1$, $\alpha=0$, $\gamma>0$), etc.

Since CKN inequalities are vital in various aspects of physics and mathematics, many scholars paid attention to them, we refer the interested readers to see \cite{Catrina09,Catrina01,Cazacu21,Cazacu22,
Cazacu23,Costa08,Do23,
Dong18,Duong21,Xia07} for CKN inequalities,
\cite{Brezis18,Cordero-Erausquin04,DelPino02,
Nguyen15} for Gagliardo-Nirenberg inequalities, \cite{Cordero-Erausquin04,DelPino03,Nguyen15,Talenti76} for Sobolev inequalities, \cite{Cazacu24,Duy22,Frank08,Lam20,Vazquez00} for Hardy inequalities, \cite{Cazacu20,Do23-2,Tertikas07} for Hardy-Rellich inequalities, etc.

\subsection{First order Caffarelli-Kohn-Nirenberg inequalities}\label{1.2}

\noindent In \cite{Xia07}, Xia studied a special case of \eqref{1.4} (with $q=p(t-1)/(p-1)$ and $\delta=1/t$), and obtained
\begin{equation}\label{1.5}
\left(\int_{\mathbb{R}^N}
\frac{|\nabla u|^p}
{|x|^{p\alpha}}\mathrm{d}x
\right)^{\frac{1}{p}}
\left[\int_{\mathbb{R}^N}
\frac{|u|^{\frac{p(t-1)}{p-1}}}
{|x|^{\frac{p(t-1)}{p-1}\beta}}
\mathrm{d}x\right]^{\frac{p-1}{p}}
\geq \frac{N-t\gamma}{t}
\int_{\mathbb{R}^N}
\frac{|u|^t}{|x|^{t\gamma}}
\mathrm{d}x,
\end{equation}
where $N\ge2$, $t>p>1$ and $\alpha,\beta,\gamma$ satisfying
$$
N-p\alpha>0,
\
N-\frac{p(t-1)}{p-1}\beta>0,
\
N-t\gamma>0
$$
and
$$
t\gamma=\alpha+(t-1)\beta+1.
$$
Moreover, when $\alpha-\frac{(t-1)}{p-1}\beta+1>0$ and
$$
N-\frac{p(t-1)}{p-1}\beta
<\left[\alpha-\frac{(t-1)\beta}{p-1}+1\right]
\frac{p(t-1)}{t-p},
$$
then $(N-t\gamma)/t$ is the optimal constant of \eqref{1.5} and can be achieved by
$$
u(x)
=\left[c+|x|^{\alpha-\frac{(t-1)\beta}{p-1}+1}
\right]^{\frac{1-p}{t-p}},
\
c>0.
$$

Especially, when $p=q=t$, $\delta=\frac{1}{p}$ and $t\gamma=\alpha+(p-1)\beta+1$ in the CKN inequality \eqref{1.4}, then \eqref{1.4} reduces into the following so-called $L^p$-CKN inequality
\begin{align}\label{1.6}
&\left(\int_{\mathbb{R}^N}
\frac{|\nabla u|^p}
{|x|^{p\alpha}}\mathrm{d}x
\right)^{\frac{1}{p}}
\left(\int_{\mathbb{R}^N}
\frac{|u|^p}{|x|^{p\beta}}
\mathrm{d}x\right)^{\frac{p-1}{p}}
\geq C(N,p,\alpha,\beta)
\int_{\mathbb{R}^N}
\frac{|u|^p}{|x|^{\alpha+(p-1)\beta+1}}
\mathrm{d}x.
\end{align}
Recently, Do \emph{et al.} \cite{Do23} studied the $L^p$-CKN inequality \eqref{1.6} and got the results as stated below.

\vspace{0.25em}
\vspace*{0.25em}

\noindent\textbf{Theorem}
(\cite[Corollary 1.2]{Do23})\textbf{.}
\emph{Assume that $N\ge1$ and $p>1$. Then there holds}
\begin{equation}\label{1.7}
\left(\int_{\mathbb{R}^N}
\frac{|\nabla u|^p}
{|x|^{p\alpha}}\mathrm{d}x
\right)^{\frac{1}{p}}
\left(\int_{\mathbb{R}^N}
\frac{|u|^p}{|x|^{p\beta}}
\mathrm{d}x\right)^{\frac{p-1}{p}}
\geq\frac{\left|N-\alpha-(p-1)\beta-1\right|}{p}
\int_{\mathbb{R}^N}
\frac{|u|^p}{|x|^{\alpha+(p-1)\beta+1}}
\mathrm{d}x,
\end{equation}
\noindent \emph{for all $u\in\mathcal{C}_0^\infty
(\mathbb{R}^N\setminus\{\mathbf{0}\})$. Moreover},
\begin{enumerate}
[itemsep=0pt, topsep=2pt, parsep=0pt]

\item[\emph{(1)}] \emph{if $\alpha-\beta+1>0$ and $\alpha\le\frac{N-p}{p}$, then $\frac{N-\alpha-(p-1)\beta-1}{p}$ is the optimal constant of \eqref{1.7} and the equality in \eqref{1.7} is only achieved by $u(x)=\Lambda e^{\frac{\lambda}{\alpha-\beta+1}
    |x|^{\alpha-\beta+1}}$, where $\Lambda\in\mathbb{R}$ and $\lambda<0$;}

\item[\emph{(2)}] \emph{if $\alpha-\beta+1<0$ and $\alpha\ge\frac{N-p}{p}$, then $\frac{\alpha+(p-1)\beta+1-N}{p}$ is the optimal constant of \eqref{1.7} and the equality in \eqref{1.7} is only achieved by $u(x)=\Lambda e^{\frac{\lambda}{\alpha-\beta+1}
    |x|^{\alpha-\beta+1}}$, where $\Lambda\in\mathbb{R}$ and $\lambda>0$.}

\end{enumerate}

When $p=2$ in \eqref{1.6}, we obtain the following $L^2$-CKN inequality
\begin{align}\label{1.8}
&\left(\int_{\mathbb{R}^N}
\frac{|\nabla u|^2}
{|x|^{2\alpha}}\mathrm{d}x
\right)^{\frac{1}{2}}
\left(\int_{\mathbb{R}^N}
\frac{|u|^2}{|x|^{2\beta}}
\mathrm{d}x\right)^{\frac{1}{2}}
\geq C(N,\alpha,\beta)
\int_{\mathbb{R}^N}
\frac{|u|^2}{|x|^{\alpha+\beta+1}}
\mathrm{d}x.
\end{align}
Recall that
\begin{equation}\label{1.9}
\begin{cases}
\mathcal{P}_1:=\left\{(\alpha,\beta):\alpha-\beta+1>0,\ \alpha\le\frac{N-2}{2}\right\};\\[0.5mm]
\mathcal{P}_2:=\left\{(\alpha,\beta):\alpha-\beta+1<0,\ \alpha\ge\frac{N-2}{2}\right\};\\[0.5mm]
\mathcal{P}:=\mathcal{P}_1\cup\mathcal{P}_2;\\[0.5mm]
\mathcal{Q}_1:=\left\{(\alpha,\beta):\alpha-\beta+1<0,\ \alpha\le\frac{N-2}{2}\right\};\\[0.5mm]
\mathcal{Q}_2:=\left\{(\alpha,\beta):\alpha-\beta+1>0,\ \alpha\ge\frac{N-2}{2}\right\};\\[0.5mm]
\mathcal{Q}:=\mathcal{Q}_1\cup\mathcal{Q}_2.
\end{cases}
\end{equation}
Costa \cite{Costa08} firstly analyzed the constant $C(N,\alpha,\beta)$ of \eqref{1.8} when  $(\alpha,\beta)\in\mathcal{P}_1$ or $(\alpha,\beta)\in\mathcal{P}_2$, and obtained that the constant $C(N,\alpha,\beta)
=\frac{|N-(\alpha+\beta+1)|}{2}$ of \eqref{1.8} is optimal. Subsequently, Catrina and Costa \cite{Catrina09} improved those results in \cite{Costa08}, the main results of \cite{Catrina09} state as follows.

\vspace{0.25em}
\vspace*{0.25em}

\noindent\textbf{Theorem}
 (\cite[Theorem 1]{Catrina09})\textbf{.}
\emph{For a nonzero constant $\Lambda$, there hold}
\begin{enumerate}
[itemsep=0pt, topsep=2pt, parsep=0pt]

\item[\emph{(1)}] \emph{if $(\alpha,\beta)\in\mathcal{P}$, $C(N,\alpha,\beta)=\frac{|N-(\alpha+\beta+1)|}{2}$ is the optimal constant of \eqref{1.8} and can be achieved by $u(x)=\Lambda e^{
    \frac{\lambda}{\alpha-\beta+1}
    |x|^{\alpha-\beta+1}},$
    where $\lambda<0$ in $\mathcal{P}_1$ and $\lambda>0$ in $\mathcal{P}_2$;}

\item[\emph{(2)}] \emph{if $(\alpha,\beta)\in\mathcal{Q}$, $C(N,\alpha,\beta)=\frac{|N-(3\alpha-\beta+3)|}{2}$ is the optimal constant of \eqref{1.8} and can be achieved by $u(x)=\Lambda|x|^{2(\alpha+1)-N}
    e^{\frac{\lambda}{\alpha-\beta+1}
    |x|^{\alpha-\beta+1}},$
    where $\lambda>0$ in $\mathcal{Q}_1$ and $\lambda<0$ in $\mathcal{Q}_2$;}

\item[\emph{(3)}] \emph{if $\beta=\alpha+1$, the optimal constant $C(N,\alpha,\alpha+1)
    =\frac{|N-2(\alpha+1)|}{2}$ cannot be achieved.}
\end{enumerate}

\vspace{0.25em}
\vspace*{0.25em}

\noindent Recently, for all $(\alpha,\beta)\in\mathbb{R}^2$ defined as \eqref{1.9}, Cazacu \emph{et al.} in \cite{Cazacu21} provided a simple proof of the following inequality to derive optimal constants and optimizers,
\begin{align}\label{1.10}
&\left(\int_{\mathbb{R}^N}
\frac{|x\cdot\nabla u|^2}
{|x|^{2\alpha+2}}\mathrm{d}x
\right)^{\frac{1}{2}}
\left(\int_{\mathbb{R}^N}
\frac{|u|^2}{|x|^{2\beta}}
\mathrm{d}x\right)^{\frac{1}{2}}
\geq \tilde{C}(N,\alpha,\beta)
\int_{\mathbb{R}^N}
\frac{|u|^2}{|x|^{\alpha+\beta+1}}
\mathrm{d}x,
\end{align}
where $\tilde{C}(N,\alpha,\beta)=C(N,\alpha,\beta)>0$ is the optimal constant of \eqref{1.10}.

In particular, the $L^2$-CKN inequality \eqref{1.8} contains HUP ($\alpha=0$, $\beta=-1$), Hydrogen Uncertainty Principles ($\alpha=\beta=0$) and Hardy inequalities ($\alpha=0$, $\beta=1$). Here we would like to emphasize that HUP play a vital role in quantum mechanics, which stated as below,
\begin{equation}\label{1.11}
\int_{\mathbb{R}^N}
|\nabla u|^2\mathrm{d}x
\int_{\mathbb{R}^N}
|x|^{2}|u|^2\mathrm{d}x
\ge\frac{N^2}{4}
\left(\int_{\mathbb{R}^N}|u|^2
\mathrm{d}x\right)^2,
\end{equation}
where $\frac{N^2}{4}$ is the optimal constant of \eqref{1.11} and can be attained by $u(x)=\Lambda e^{-\lambda|x|^2}$ for $\Lambda\in\mathbb{R}$ and $\lambda>0$, see \cite{Folland97,Lieb10,Weyl50} for more details.

\subsection{Second order Caffarelli-Kohn-Nirenberg inequalities}\label{1.3}

\noindent Compared with the first order CKN inequalities, there are few results (see, e.g. \cite{Cazacu22,Cazacu23,Duong21,Lin86}) about the higher order CKN inequalities. In \cite{Cazacu22}, Cazacu \emph{et al.} got the following second order HUP, for $N\ge1$ and $u\in\mathcal{S}(\mathbb{R}^N)$ (which is the Schwartz space of smooth functions),
\begin{equation}\label{1.12}
\int_{\mathbb{R}^N}
|\Delta u|^2\mathrm{d}x
\int_{\mathbb{R}^N}
|x|^{2}|\nabla u|^2\mathrm{d}x
\ge\frac{\left(N+2\right)^2}{4}
\left(\int_{\mathbb{R}^N}|\nabla u|^2
\mathrm{d}x\right)^2,
\end{equation}
where the constant $\frac{(N+2)^2}{4}$ of \eqref{1.12} is optimal and can be attained by $u(x)=\Lambda e^{-\lambda|x|^2}$ for $\Lambda\in\mathbb{C}$ and $\lambda>0$. Later, in \cite{Duong21}, Duong and Nguyen improved the result of \cite{Cazacu22} to the following weighted version,
\begin{equation}\label{1.13}
\int_{\mathbb{R}^N}
\frac{|\Delta u|^2}{|x|^{2\alpha}}\mathrm{d}x
\int_{\mathbb{R}^N}
|x|^{2\alpha}|\nabla u\cdot x|^2\mathrm{d}x
\ge\frac{\left(N+4\alpha+2\right)^2}{4}
\left(\int_{\mathbb{R}^N}|\nabla u|^2
\mathrm{d}x\right)^2,
\end{equation}
where $N\ge1$, $N-2\alpha>0$, $N+2\alpha>0$ and $N+4\alpha+2>0$. Moreover, if $\alpha+1>0$, $\frac{(N+4\alpha+2)^2}{4}$ is the optimal constant of \eqref{1.13} and there is $u(x)=\exp\left[-\frac{1}
{2(\alpha+1)}|x|^{2(\alpha+1)}\right]$ such that the optimal constant can be attained. The authors in  \cite{Duong21} also obtained the following inequality, for all radial functions $u\in\mathcal{C}_0^\infty(\mathbb{R}^N
\setminus\{\mathbf{0}\})$,
\begin{equation}\label{1.14}
\left(\int_{\mathbb{R}^N}
\frac{|\Delta u|^2}
{|x|^{2\alpha}}\mathrm{d}x
\right)^{\frac{1}{2}}
\left[\int_{\mathbb{R}^N}
\frac{|\nabla u|^{2(t-1)}}
{|x|^{2(t-1)\beta}}
\mathrm{d}x\right]^{\frac{1}{2}}
\geq \frac{\left|N+t\left(2\alpha
-\gamma+1\right)\right|}{t}
\int_{\mathbb{R}^N}
\frac{|\nabla u|^t}{|x|^{t\gamma}}
\mathrm{d}x,
\end{equation}
where $N-2\alpha>0$, $N-2(t-1)\beta>0$, $N-t\gamma>0$ and $t\gamma=\alpha+(t-1)\beta+1$. Provided that some suitable assumptions, they got that $\frac{N+t\left(2\alpha
-\gamma+1\right)}{t}$ is the optimal constant of \eqref{1.14}, and if $t=2$, the optimal constant is attained only by
$$
u(x)=\int_{|x|}^\infty
r^{2\alpha+1}\exp\left[{-\frac{r^{\alpha-(t-1)\beta+1}}
{\alpha-(t-1)\beta+1}}\right]\mathrm{d}r;
$$
if $t>2$, the equality in \eqref{1.14} holds if and only if
$$
u(x)=\int_{|x|}^\infty
r^{2\alpha+1}\left[1+
\frac{(t-2)
r^{\alpha-(t-1)\beta+1+(2\alpha+1)(t-2)}}
{\alpha-(t-1)\beta+1+(2\alpha+1)(t-2)
}\right]^{\frac{1}{2-t}}\mathrm{d}r.
$$

Motivated by the above results, it is natural to ask some questions.

\vspace{0.5em}

\textbf{$\bullet$\,\emph{Q1}}: We are curious whether  inequality \eqref{1.5} has a second order version. If this inequality exists, are there optimal constants and optimizers of it.

\textbf{$\bullet$\,\emph{Q2:}} We wonder if there exists a second order version of inequality (1.7). If this inequality exists, and do there exist optimal constants and optimizers of it?

\textbf{$\bullet$\,\emph{Q3:}} Under some more general assumptions, do there exist optimal constants and optimizers of \eqref{1.13}?

\vspace{0.5em}

We give affirmative answers to Q1, Q2 and Q3 in Theorems \ref{thm-2.1}, \ref{thm-2.5} and \ref{thm-2.8}, respectively.

\section{Main results}

\noindent Main results and some relevant comments will be presented in this section. Let $N\ge1$, $p,q\ge1$ and
$$
\mathcal{H}_{\alpha,\beta}^{p,q}(\mathbb{R}^N):
=\left\{u\in\mathcal{C}_0^\infty
\left(\mathbb{R}^N\setminus\{\mathbf{0}\}\right):
\int_{\mathbb{R}^N}\frac{|\Delta u|^p}
{|x|^{p\alpha}}\mathrm{d}x<\infty,
\
\int_{\mathbb{R}^N}\frac{|\nabla u|^q}
{|x|^{q\beta}}\mathrm{d}x<\infty
\right\}.
$$

\subsection{Optimal constants and optimizers for a more general subfamily of the second order CKN inequalities}

\noindent We firstly consider the following inequality \eqref{2.2} for radial functions.

\begin{thm}\label{thm-2.1}
Assume that $N\ge1$, $t\ge p>1$ and $\alpha,\beta,\gamma$ satisfy
\begin{equation}\label{2.1}
t\gamma=\alpha+(t-1)\beta+1,
\end{equation}
then for all radial functions
$u\in\mathcal{H}_{\alpha,\beta}^{p,\frac{p(t-1)}{p-1}}
(\mathbb{R}^N)$,
\begin{equation}\label{2.2}
\left(\int_{\mathbb{R}^N}
\frac{|\Delta u|^p}{|x|^{p\alpha}}
\mathrm{d}x\right)^{\frac{1}{p}}
\left[\int_{\mathbb{R}^N}
\frac{\left|\nabla u\right|^{\frac{p(t-1)}{p-1}}}
{|x|^{\frac{p(t-1)}{p-1}\beta}}
\mathrm{d}x\right]^{\frac{p-1}{p}}
\ge\frac{\left|N-t(N+\gamma-1)\right|}{t}
\int_{\mathbb{R}^N}
\frac{\left|\nabla u\right|^t}{|x|^{t\gamma}}
\mathrm{d}x.
\end{equation}
Moreover,
\begin{itemize}
\item Case 1: $p=t$.
\begin{enumerate}[itemsep=0pt, topsep=0pt, parsep=0pt]

\item[(1a)] If $\alpha-\beta+1>0$ and $p\alpha+(p-1)N<0$, the constant $\frac{|N-t(N+\gamma-1)|}{t}
    =\frac{-\alpha-(p-1)\beta-(p-1)(N-1)}{p}$ is optimal, and the equality in \eqref{2.2} holds if and only if
$$
U_1(x)=\Lambda\int_{|x|}^\infty
r^{1-N}\exp\left(-\frac{\lambda r^{\alpha-\beta+1}}
{\alpha-\beta+1}\right)
\mathrm{d}r,
$$
for $\Lambda\in\mathbb{R}$ and $\lambda>0$;

\item[(1b)] if $\alpha-\beta+1<0$ and $p\alpha+(p-1)N>0$, the constant
    $\frac{|N-t(N+\gamma-1)|}{t}
    =\frac{\alpha+(p-1)\beta+(p-1)(N-1)}{p}$ is optimal, and the equality in \eqref{2.2} holds if and only if

$$
U_2(x)=\Lambda\int_{|x|}^\infty
r^{1-N}\exp\left(\frac{\lambda r^{\alpha-\beta+1}}
{\alpha-\beta+1}\right)
\mathrm{d}r,
$$
for $\Lambda\in\mathbb{R}$ and $\lambda>0$.
\end{enumerate}

\item Case 2: $p<t$.
Under the following assumptions
\begin{align*}
t-p&<(p-1)(\alpha+1)-(t-1)\beta,\\
N-p\alpha
&<\frac{p\left[(p-1)\alpha-(t-1)\beta+t-1\right]}{t-p},
\\
N-\frac{p(t-1)}{p-1}\beta
&<\frac{p(t-1)
\left[(p-1)(\alpha+1)-(t-1)\beta\right]}{(t-p)(p-1)},
\end{align*}
\begin{enumerate}[itemsep=0pt, topsep=0pt, parsep=0pt]

\item[(2a)] if $(p-1)(\alpha+1)-(t-1)\beta-(N-1)(t-p)>0$ and
    $p\alpha+(p-1)N<0$, the constant $\frac{|N-t(N+\gamma-1)|}{t}
    =\frac{N-t(N+\gamma-1)}{t}$ is optimal, and the equality in \eqref{2.2} holds if and only if
$$
U_3(x)=\Lambda\int_{|x|}^\infty
r^{1-N}
\left[1+
\frac{(t-p)\lambda r^{\frac{(p-1)(\alpha+1)-(t-1)\beta-(N-1)(t-p)}{p-1}}}
{(p-1)(\alpha+1)-(t-1)\beta-(N-1)(t-p)}
\right]^{\frac{p-1}{p-t}}
\mathrm{d}r,
$$
for $\Lambda\in\mathbb{R}$ and $\lambda>0$;

\item[(2b)] if $(p-1)(\alpha+1)-(t-1)\beta-(N-1)(t-p)<0$ and
    $p\alpha+(p-1)N>0$, the constant $\frac{|N-t(N+\gamma-1)|}{t}
    =\frac{t(N+\gamma-1)-N}{t}$ is optimal, and the equality in \eqref{2.2} holds if and only if
$$
U_4(x)=\Lambda\int_{|x|}^\infty
r^{1-N}
\left[1+
\frac{(p-t)\lambda r^{\frac{(p-1)(\alpha+1)-(t-1)\beta-(N-1)(t-p)}{p-1}}}
{(p-1)(\alpha+1)-(t-1)\beta-(N-1)(t-p)}
\right]^{\frac{p-1}{p-t}}
\mathrm{d}r,
$$
for $\Lambda\in\mathbb{R}$ and $\lambda>0$.

\end{enumerate}
\end{itemize}
\end{thm}

\begin{rem}
Here we present some comments about Theorem \ref{thm-2.1}.
\begin{enumerate}[itemsep=0pt, topsep=0pt, parsep=0pt]

\item[(1)] Compared with \cite{Cazacu22,Duong21}, we study a more general class of the second order CKN inequalities (with $p\neq t$ and $p\neq 2$).

\item[(2)] We also obtain the sufficient and necessary conditions such that the equality in \eqref{2.2} holds when $p=t$ and $p<t$. In order to obtain the sufficient and necessary conditions, we apply H\"{o}lder inequality.

\end{enumerate}
\end{rem}

\subsection{Second order $L^p$-CKN identities, optimal constants and optimizers for the second order $L^p$-CKN inequalities}

\noindent Now, we establish second order $L^p$-CKN identities and the optimal second order $L^p$-CKN inequalities for any functions without radial assumption.

Inspired by \cite{Do23,Duy22}, let $p>1$, $\overrightarrow{X},\overrightarrow{Y}$ be vectors in $\mathbb{R}^N$ ($N\ge1$), and
\begin{equation}\label{2.3}
\mathcal{K}_p\left(\overrightarrow{X},
\overrightarrow{Y}\right)
:=\left|\overrightarrow{Y}\right|^p
+(p-1)\left|\overrightarrow{X}\right|^p
-p\left|\overrightarrow{X}\right|^{p-2}\overrightarrow{X}
\cdot\overrightarrow{Y}.
\end{equation}

\begin{thm}\label{thm-2.3}
Let $N\ge1$, $p>1$, $\alpha-\beta+1>0$ and $p\alpha+(p-1)N<0$. For all $u\in \mathcal{H}_{\alpha,\beta}^{p,p}(\mathbb{R}^N)$, there hold
\begin{align}\label{2.4}
&\int_{\mathbb{R}^N}\frac{|\Delta u|^p}
{|x|^{p\alpha}}\mathrm{d}x
+(p-1)\int_{\mathbb{R}^N}\frac{|\nabla u|^p}
{|x|^{p\beta}}\mathrm{d}x
-\left[-\alpha-(p-1)\beta-(p-1)(N-1)\right]
\int_{\mathbb{R}^N}
\frac{|\nabla u|^p}{|x|^{\alpha+(p-1)\beta+1}}\mathrm{d}x
\nonumber\\&\qquad
+p\left[\alpha+(p-1)\beta+1\right]
\int_{\mathbb{R}^N}
\frac{|\nabla u|^{p-2}\left[\left(x\cdot\nabla u\right)^2-|x|^2|\nabla u|^2\right]}{|x|^{\alpha+(p-1)\beta+3}}\mathrm{d}x
\nonumber\\&\quad
=\int_{\mathbb{R}^N}\frac{1}
{|x|^{p\alpha}}\mathcal{K}_p
\left(-|x|^{\alpha-\beta-1}x\cdot\nabla u,\Delta u\right)\mathrm{d}x,
\end{align}
and
\begin{align}\label{2.5}
&\left(\int_{\mathbb{R}^N}\frac{|\Delta u|^p}
{|x|^{p\alpha}}\mathrm{d}x\right)^{\frac{1}{p}}
\left(\int_{\mathbb{R}^N}\frac{|\nabla u|^p}
{|x|^{p\beta}}\mathrm{d}x\right)^{\frac{p-1}{p}}
\!-\!\frac{-\alpha\!-\!(p-1)\beta\!-\!(p-1)(N-1)}{p}
\int_{\mathbb{R}^N}
\frac{|\nabla u|^p}
{|x|^{\alpha+(p-1)\beta+1}}\mathrm{d}x
\nonumber\\&\qquad
+\left[\alpha+(p-1)\beta+1\right]
\int_{\mathbb{R}^N}
\frac{|\nabla u|^{p-2}\left[\left(x\cdot\nabla u\right)^2-|x|^2|\nabla u|^2\right]}{|x|^{\alpha+(p-1)\beta+3}}\mathrm{d}x
\nonumber\\&\quad
=\frac{1}{p}\int_{\mathbb{R}^N}\frac{1}
{|x|^{p\alpha}}\mathcal{K}_p
\left(-\left(\frac{\int_{\mathbb{R}^N}\frac{|\Delta u|^p}
{|x|^{p\alpha}}\mathrm{d}x}
{\int_{\mathbb{R}^N}\frac{|\nabla u|^p}
{|x|^{p\beta}}\mathrm{d}x}\right)^{\frac{1}{p^2}}
|x|^{\alpha-\beta-1}x\cdot\nabla u,\left(\frac{\int_{\mathbb{R}^N}\frac{|\nabla u|^p}
{|x|^{p\beta}}\mathrm{d}x}
{\int_{\mathbb{R}^N}\frac{|\Delta u|^p}
{|x|^{p\alpha}}\mathrm{d}x}\right)^{\frac{p-1}{p^2}}\Delta u\right)\mathrm{d}x.
\end{align}
\end{thm}

\begin{thm}\label{thm-2.4}
Let $N\ge1$, $p>1$, $\alpha-\beta+1<0$ and $p\alpha+(p-1)N>0$. Then for all $u\in \mathcal{H}_{\alpha,\beta}^{p,p}(\mathbb{R}^N)$,
\begin{align}\label{2.6}
&\int_{\mathbb{R}^N}\frac{|\Delta u|^p}
{|x|^{p\alpha}}\mathrm{d}x
+(p-1)\int_{\mathbb{R}^N}\frac{|\nabla u|^p}
{|x|^{p\beta}}\mathrm{d}x
-\left[\alpha+(p-1)\beta+(p-1)(N-1)\right]
\int_{\mathbb{R}^N}
\frac{|\nabla u|^p}{|x|^{\alpha+(p-1)\beta+1}}\mathrm{d}x
\nonumber\\&\qquad
-p\left[\alpha+(p-1)\beta+1\right]
\int_{\mathbb{R}^N}
\frac{|\nabla u|^{p-2}\left[\left(x\cdot\nabla u\right)^2-|x|^2|\nabla u|^2\right]}{|x|^{\alpha+(p-1)\beta+3}}\mathrm{d}x
\nonumber\\&\quad
=\int_{\mathbb{R}^N}\frac{1}
{|x|^{p\alpha}}\mathcal{K}_p
\left(|x|^{\alpha-\beta-1}x\cdot\nabla u,\Delta u\right)\mathrm{d}x,
\end{align}
and
\begin{align}\label{2.7}
&\left(\int_{\mathbb{R}^N}\frac{|\Delta u|^p}
{|x|^{p\alpha}}\mathrm{d}x\right)^{\frac{1}{p}}
\left(\int_{\mathbb{R}^N}\frac{|\nabla u|^p}
{|x|^{p\beta}}\mathrm{d}x\right)^{\frac{p-1}{p}}
-\frac{\alpha+(p-1)\beta+(p-1)(N-1)}{p}
\int_{\mathbb{R}^N}
\frac{|\nabla u|^p}{|x|^{\alpha+(p-1)\beta+1}}\mathrm{d}x
\nonumber\\&\qquad
-\left[\alpha+(p-1)\beta+1\right]
\int_{\mathbb{R}^N}
\frac{|\nabla u|^{p-2}\left[\left(x\cdot\nabla u\right)^2-|x|^2|\nabla u|^2\right]}{|x|^{\alpha+(p-1)\beta+3}}\mathrm{d}x
\nonumber\\&\quad
=\frac{1}{p}\int_{\mathbb{R}^N}\frac{1}
{|x|^{p\alpha}}\mathcal{K}_p
\left(\left(\frac{\int_{\mathbb{R}^N}\frac{|\Delta u|^p}
{|x|^{p\alpha}}\mathrm{d}x}
{\int_{\mathbb{R}^N}\frac{|\nabla u|^p}
{|x|^{p\beta}}\mathrm{d}x}\right)^{\frac{1}{p^2}}
|x|^{\alpha-\beta-1}x\cdot\nabla u,\left(\frac{\int_{\mathbb{R}^N}\frac{|\nabla u|^p}
{|x|^{p\beta}}\mathrm{d}x}
{\int_{\mathbb{R}^N}\frac{|\Delta u|^p}
{|x|^{p\alpha}}\mathrm{d}x}\right)^{\frac{p-1}{p^2}}\Delta u\right)\mathrm{d}x.
\end{align}
\end{thm}

\begin{thm}\label{thm-2.5}
Assume that $N\ge1$ and $p>1$. For all $u\in \mathcal{H}_{\alpha,\beta}^{p,p}(\mathbb{R}^N)$,
\begin{enumerate}[itemsep=0pt, topsep=0pt, parsep=0pt]

\item[(1)] if $\alpha+(p-1)\beta+1>0$, $\alpha-\beta+1>0$ and $p\alpha+(p-1)N<0$, then
\begin{eqnarray}\label{2.8}
\!\left(\!\int_{\mathbb{R}^N}\!\!\frac{|\Delta u|^p}
{|x|^{p\alpha}}\mathrm{d}x\!\!\right)^{\frac{1}{p}}
\!\!\left(\!\int_{\mathbb{R}^N}\!\!\frac{|\nabla u|^p}
{|x|^{p\beta}}\mathrm{d}x\!\!\right)^{\frac{p-1}{p}}\!
\!\!
\ge\!\frac{-\alpha\!-\!(p\!-\!1)\beta
\!-\!(p\!-\!1)(N\!\!-\!1)}{p}
\!\!\int_{\mathbb{R}^N}\!\!
\frac{|\nabla u|^p}{|x|^{\alpha+(p-\!1)\beta+1}}\mathrm{d}x,
\end{eqnarray}
the constant $\frac{-\alpha-(p-1)\beta-(p-1)(N-1)}{p}$ is optimal, and \eqref{2.8} is reached with equality if and only if
$$
U_5(x)=\Lambda\int_{|x|}^\infty
r^{1-N}\exp\left(-\frac{\lambda r^{\alpha-\beta+1}}
{\alpha-\beta+1}\right)
\mathrm{d}r,
$$
with $\Lambda\in\mathbb{R}$ and $\lambda>0$;

\item[(2)]
if $\alpha+(p-1)\beta+1<0$,
$\alpha-\beta+1<0$ and
$p\alpha+(p-1)N>0$, then
\begin{eqnarray}\label{2.9}
\!\left(\!\int_{\mathbb{R}^N}\!\!\frac{|\Delta u|^p}
{|x|^{p\alpha}}\mathrm{d}x\!\!\right)^{\frac{1}{p}}
\!\!\left(\!\int_{\mathbb{R}^N}\!\!\frac{|\nabla u|^p}
{|x|^{p\beta}}\mathrm{d}x\!\!\right)^{\frac{p-1}{p}}\!
\!\!
\ge\!\frac{\alpha\!+\!(p\!-\!1)\beta
\!+\!(p\!-\!1)(N\!\!-\!1)}{p}
\!\!\int_{\mathbb{R}^N}\!\!
\frac{|\nabla u|^p}{|x|^{\alpha+(p-\!1)\beta+1}}\mathrm{d}x,
\end{eqnarray}
the constant $\frac{\alpha+(p-1)\beta+(p-1)(N-1)}{p}$ is optimal, and \eqref{2.9} is reached with equality if and only if
$$
U_6(x)=\Lambda\int_{|x|}^\infty
r^{1-N}\exp\left(\frac{\lambda r^{\alpha-\beta+1}}
{\alpha-\beta+1}\right)
\mathrm{d}r,
$$
with $\Lambda\in\mathbb{R}$ and $\lambda>0$;

\item[(3)] if $\alpha+(p-1)\beta+1=0$, there holds
\begin{eqnarray}\label{2.14}
\!\left(\!\int_{\mathbb{R}^N}\!\!\frac{|\Delta u|^p}
{|x|^{p\alpha}}\mathrm{d}x\!\!\right)^{\frac{1}{p}}
\!\!\left(\!\int_{\mathbb{R}^N}\!\!\frac{|\nabla u|^p}
{|x|^{p\beta}}\mathrm{d}x\!\!\right)^{\frac{p-1}{p}}\!
\!\!
\ge\!\frac{\left|\alpha\!+\!(p\!-\!1)\beta
\!+\!(p\!-\!1)(N\!\!-\!1)\right|}{p}
\!\!\int_{\mathbb{R}^N}\!\!
\frac{|\nabla u|^p}{|x|^{\alpha+(p-\!1)\beta+1}}\mathrm{d}x.
\end{eqnarray}

\noindent Moreover,

\begin{enumerate}[itemsep=0pt, topsep=0pt, parsep=0pt]

\item[(3a)] if $\alpha-\beta+1>0$ and $p\alpha+(p-1)N<0$, then $\frac{-\alpha-(p-1)\beta-(p-1)(N-1)}{p}$ is the optimal constant of \eqref{2.14} and there are nontrivial functions such that the optimal constant can be attained;

\item[(3b)] if $\alpha-\beta+1<0$ and $p\alpha+(p-1)N>0$, then
    $\frac{\alpha+(p-1)\beta+(p-1)(N-1)}{p}$ is the  optimal constant of \eqref{2.14} and there are nontrivial functions such that the optimal constant can be attained.

\end{enumerate}

\end{enumerate}
\end{thm}

\begin{thm}\label{thm-2.6}
Assume that $N\ge1$ and $p>1$. For all $u\in \mathcal{H}_{\alpha,\beta}^{p,p}(\mathbb{R}^N)$,
\begin{enumerate}[itemsep=0pt, topsep=0pt, parsep=0pt]

\item[(1)] if $\alpha+(p-1)\beta+1>0$, $\alpha-\beta+1>0$ and $p\alpha+(p-1)N<0$, then
\begin{eqnarray}\label{2.10}
\!\!\!\int_{\mathbb{R}^N}\!\!\frac{|\Delta u|^p}
{|x|^{p\alpha}}\mathrm{d}x
\!+\!(p\!-\!1)\!\!\int_{\mathbb{R}^N}\!\!\frac{|\nabla u|^p}{|x|^{p\beta}}\mathrm{d}x
\!\ge\!\left[-\alpha\!
-\!(p\!-\!1)\beta\!\!-\!(p\!-\!1)(N\!\!-\!1)\right]\!
\!\int_{\mathbb{R}^N}\!\!
\frac{|\nabla u|^p}{|x|^{\alpha+(p-1)\beta+1}}\mathrm{d}x,
\end{eqnarray}
the constant $-\alpha-(p-1)\beta-(p-1)(N-1)$ is optimal, and \eqref{2.10} is reached with equality if and only if
$$
U_7(x)=\Lambda\int_{|x|}^\infty
r^{1-N}\exp\left(-\frac{r^{\alpha-\beta+1}}
{\alpha-\beta+1}\right)
\mathrm{d}r,
$$
for $\Lambda\in\mathbb{R}$;

\item[(2)] if $\alpha+(p-1)\beta+1<0$, $\alpha-\beta+1<0$ and $p\alpha+(p-1)N>0$, then
\begin{eqnarray}\label{2.11}
\!\!\int_{\mathbb{R}^N}\!\!\frac{|\Delta u|^p}
{|x|^{p\alpha}}\mathrm{d}x
\!+\!(p\!-\!1)\!\!\int_{\mathbb{R}^N}\!\!\frac{|\nabla u|^p}{|x|^{p\beta}}\mathrm{d}x
\!\ge\!\left[\alpha\!
+\!(p\!-\!1)\beta\!+\!(p\!-\!1)(N\!\!-\!1)\right]\!
\!\int_{\mathbb{R}^N}\!\!
\frac{|\nabla u|^p}{|x|^{\alpha+(p-1)\beta+1}}\mathrm{d}x,
\end{eqnarray}
the constant $\alpha+(p-1)\beta+(p-1)(N-1)$ is optimal, and \eqref{2.11} is reached with equality if and only if
$$
U_8(x)=\Lambda\int_{|x|}^\infty
r^{1-N}\exp\left(\frac{r^{\alpha-\beta+1}}
{\alpha-\beta+1}\right)
\mathrm{d}r,
$$
for $\Lambda\in\mathbb{R}$;

\item[(3)] if $\alpha+(p-1)\beta+1=0$, there holds
\begin{eqnarray}\label{2.15}
\!\!\int_{\mathbb{R}^N}\!\!\frac{|\Delta u|^p}
{|x|^{p\alpha}}\mathrm{d}x
\!+\!(p\!-\!1)\!\!\int_{\mathbb{R}^N}\!\!\frac{|\nabla u|^p}{|x|^{p\beta}}\mathrm{d}x
\!\ge\!\left|\alpha\!
+\!(p\!-\!1)\beta\!+\!(p\!-\!1)(N\!\!-\!1)\right|\!
\!\int_{\mathbb{R}^N}\!\!
\frac{|\nabla u|^p}{|x|^{\alpha+(p-1)\beta+1}}\mathrm{d}x.
\end{eqnarray}

\noindent Moreover,

\begin{enumerate}[itemsep=0pt, topsep=0pt, parsep=0pt]

\item[(3a)] if $\alpha-\beta+1>0$ and $p\alpha+(p-1)N<0$, then   $-\alpha-(p-1)\beta-(p-1)(N-1)$ is the optimal constant of \eqref{2.15} and there exist some nontrivial functions such that the optimal constant can be reached;

\item[(3b)] if $\alpha-\beta+1<0$ and $p\alpha+(p-1)N>0$, then $\alpha+(p-1)\beta+(p-1)(N-1)$ is the optimal constant of \eqref{2.15} and there exist some nontrivial functions such that the optimal constant can be reached.

\end{enumerate}

\end{enumerate}
\end{thm}

\begin{rem}
Some comments about Theorems \ref{thm-2.3}$-$\ref{thm-2.6} state in the following three aspects.
\begin{enumerate}[itemsep=0pt, topsep=0pt, parsep=0pt]

\item[(1)] We establish the second order $L^p$-CKN identities, which can be used to study the stability estimates of the second order $L^p$-CKN inequalities in the follow-up study.

\item[(2)] Compared with Theorem \ref{thm-2.1} when $p=t$, by applying the second order $L^p$-CKN identities and the negativity of $\mathcal{K}_p$, we obtain the same results without radial assumption.

\item[(3)] The sufficient and necessary conditions are shown to guarantee that the equality in \eqref{2.8}-\eqref{2.15} hold. This is our main difficulty. Based on some properties of $\mathcal{K}_p$ and the condition such that the equality of $(x\cdot\nabla u)^2\le|x|^2|\nabla u|^2$ holds, we overcome it.

\end{enumerate}
\end{rem}

\subsection{Optimal constants and optimizers for the weighted second order HUP}

\begin{thm}\label{thm-2.8}
Let $N\ge1$ and $\alpha\in\mathbb{R}$ satisfy one of the following three assumptions:
\begin{enumerate}[itemsep=0pt, topsep=0pt, parsep=0pt]

\setlength{\itemindent}
{2em}\item[Case 1:] $\alpha+1>0$, $N-2\alpha>0$
and $N-(\sqrt{5}-1)\alpha-(\sqrt{5}+1)>0$;

\setlength{\itemindent}
{2em}\item[Case 2:] $N+2\alpha<0$,  $N+(\sqrt{5}+1)\alpha+(\sqrt{5}-1)>0$ and $N-2\sqrt{2}\alpha-2\sqrt{2}>0$;

\setlength{\itemindent}
{2em}\item[Case 3:] $\alpha+1<0$ and $N+4\alpha+2>0$,

\end{enumerate}
there holds
\begin{equation}\label{2.12}
\left(\int_{\mathbb{R}^N}
\frac{|\Delta u|^2}{|x|^{2\alpha}}\mathrm{d}x
\right)^{\frac{1}{2}}
\left(\int_{\mathbb{R}^N}
|x|^{2\alpha+2}|\nabla u|^2\mathrm{d}x
\right)^{\frac{1}{2}}
\ge\frac{|N+4\alpha+2|}{2}
\int_{\mathbb{R}^N}|\nabla u|^2\mathrm{d}x,
\end{equation}
for all $u\in \mathcal{H}_{\alpha,-\alpha-1}^{2,2}(\mathbb{R}^N)$. Furthermore,

\begin{enumerate}[itemsep=0pt, topsep=0pt, parsep=0pt]

\item[(1)] if $\alpha+1>0$ and $N+2\alpha>0$, $\frac{N+4\alpha+2}{2}$ is the optimal constant of \eqref{2.12} and can be achieved by some nontrivial functions, such as
    $$
    u(x)=e^{-\frac{|x|^{2(\alpha+1)}}{2(\alpha+1)}};
    $$

\item[(2)] if $\alpha+1<0$ and $N+2\alpha<0$,  $-\frac{N+4\alpha+2}{2}$ is the optimal constant of \eqref{2.12} and can be achieved by some nontrivial functions, such as
    $$
    u(x)=e^{\frac{|x|^{2(\alpha+1)}}{2(\alpha+1)}}.
    $$

\end{enumerate}
\end{thm}

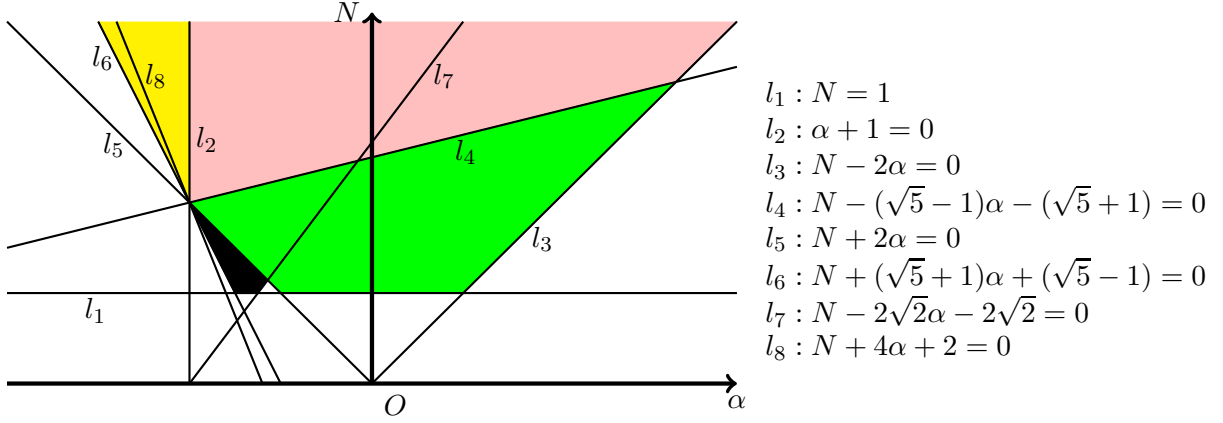
\begin{figure}
\centering
\begin{tikzpicture}[scale=2.4]

\draw[fill=pink](-1,2)--(-1,1)--(5/3,5/3)
--(2,2);
\draw[fill=yellow](-3/2,2)--(-1,1)--(-1,2);
\draw[fill=black](-1,1)--(-3/4,1/2)--(-5/8,1/2)
--(-4/7,4/7);
\draw[fill=green](-1,1)--(-1/2,1/2)--(1/2,1/2)
--(5/3,5/3);


\draw[-,thick](-2,1/2)--(2,1/2);
{\node[above left] at(-1.4,0.28){$l_1$};}

\draw[-,thick](-1,0)--(-1,2);
{\node[right] at(-1.02,1.35){$l_2$};}

\draw[-,thick](0,0)--(2,2);
{\node[left] at(1.05,0.8){$l_3$};}

\draw[-,thick](-2,3/4)--(2,7/4);
{\node[right] at(0.4,1.28){$l_4$};}

\draw[-,thick](-2,2)--(0,0);
{\node[left] at(-1.3,1.32){$l_5$};}

\draw[-,thick](-3/2,2)--(-1/2,0);
{\node[left] at(-1.36,1.8){$l_6$};}

\draw[-,thick](-1,0)--(1/2,2);
{\node[right] at(0.28,1.7){$l_7$};}

\draw[-,thick](-7/5,2)--(-3/5,0);
{\node[right] at(-1.3,1.7){$l_8$};}

		\draw[->,ultra thick](-2,0)--(0,0)node[below right]{$O$}--(2,0) node[below]{$\alpha$};
		\draw[->,ultra thick](0,0)--(0,2.05)node[left]{$N$};

\node[right] at (2.1,1.6){$l_1: N=1$};
\node[right] at (2.1,1.4){$l_2: \alpha+1=0$};
\node[right] at (2.1,1.2){$l_3: N-2\alpha=0$};
\node[right] at (2.1,1){$l_4: N-(\sqrt{5}-1)\alpha-(\sqrt{5}+1)=0$};
\node[right] at (2.1,0.8){$l_5: N+2\alpha=0$};
\node[right] at (2.1,0.6){$l_6: N+(\sqrt{5}+1)\alpha+(\sqrt{5}-1)=0$};
\node[right] at (2.1,0.4){$l_7: N-2\sqrt{2}\alpha-2\sqrt{2}=0$};
\node[right] at (2.1,0.2){$l_8: N+4\alpha+2=0$};

\node[right] at (-1.4,-0.4){\emph{Case 1}: pink};
\node[right] at (-1.4,-0.6){\emph{Case 2}: black};
\node[right] at (-1.4,-0.8){\emph{Case 3}: yellow};
\node[right] at (-1.4,-1){$``N\geq 1$, $N-2\alpha>0$, $N+2\alpha>0$, $N+4\alpha+2>0$'': pink+yellow+green};

\end{tikzpicture}
\caption{\small Choice of $(\alpha,N)$}
\label{figure1}
\end{figure}

\begin{rem}
Some comments about Theorem \ref{thm-2.8} are stated in the following.
\begin{enumerate}[itemsep=0pt, topsep=0pt, parsep=0pt]

\item[(1)] In \cite[Theorem 1.1]{Duong21}, the authors obtained \eqref{2.12} when $N\ge1$ and $\alpha\in\mathbb{R}$ satisfy ``$N-2\alpha>0$, $N+2\alpha>0$ and $N+4\alpha+2>0$''. As shown in  Figure \ref{figure1}, we see that Case 2 (the part marked in black) does not include in the region ``$N-2\alpha>0$, $N+2\alpha>0$ and $N+4\alpha+2>0$'' (the area marked in pink, yellow and green), this means that our results can be regarded as supplementary results of \cite{Duong21}.

\item[(2)] We also analyzed the optimizers of \eqref{2.12} when $\alpha+1>0$ and $N+2\alpha>0$, and $\alpha+1<0$ and $N+2\alpha<0$. It is natural to ask an open question about the existence of the optimizers of \eqref{2.12} when $\alpha+1>0$ and $N+2\alpha<0$, or $\alpha+1<0$ and $N+2\alpha>0$.

\item[(3)] The main difficulty in this result is to compute the optimal constant $\frac{|N+4\alpha+2|}{2}$ of \eqref{2.12}. We resolve it by means of some monotone properties, fine calculation and analysis (see Lemma \ref{lem-3.2} below).

\end{enumerate}
\end{rem}

\subsection{Structure of the rest of this paper}

\noindent In Sect. \ref{3.1}, we aim to consider optimal constants and optimizers for a more general subfamily of the second order CKN inequalities (proof of Theorem \ref{thm-2.1}). The second order $L^p$-CKN identities and inequalities are established in Sect. \ref{3.2} (proof of Theorems \ref{thm-2.3}, \ref{thm-2.4}, \ref{thm-2.5} and \ref{thm-2.6}). In Sect. \ref{3.3}, we study the optimal weighted second order HUP under more general assumptions (proof of Theorem \ref{thm-2.8}).

\section{Proof of main results}\label{3}

\subsection{Optimal constants and optimizers for a more general subfamily of the second order CKN inequalities: proof of Theorem \ref{thm-2.1}}\label{3.1}

\begin{proof}[\rm\textbf{Proof of Theorem \ref{thm-2.1}}]
The density argument allows us to prove \eqref{2.2} for radial functions $u\in\mathcal{C}_0^\infty(\mathbb{R}^N
\setminus\{\mathbf{0}\})$. By using integration by parts,
\begin{align*}
\int_{\mathbb{R}^N}
\frac{\left|\nabla u\right|^t}{|x|^{t\gamma}}
\mathrm{d}x
&=
\left|\mathbb{S}^{N-1}\right|
\int_0^\infty\left|u'\right|^t
r^{N-t\gamma-1}\mathrm{d}r
\\&=
\frac{1}{N-t\gamma}\left|\mathbb{S}^{N-1}\right|
\int_0^\infty\left|u'\right|^t
\left(r^{N-t\gamma}\right)'\mathrm{d}r
\\&=
-\frac{t}{N-t\gamma}\left|\mathbb{S}^{N-1}\right|
\int_0^\infty|u'|^{t-2}u'u''r^{N-t\gamma}\mathrm{d}r
\\&=
-\frac{t}{N-t\gamma}\left|\mathbb{S}^{N-1}\right|
\int_0^\infty\left|u'\right|^{t-2}u'
\left(u''+\frac{N-1}{r}u'\right)
r^{N-t\gamma}\mathrm{d}r
\\&\quad+
\frac{t(N-1)}{N-t\gamma}\left|\mathbb{S}^{N-1}\right|
\int_0^\infty\left|u'\right|^tr^{N-t\gamma-1}\mathrm{d}r
\\&=
-\frac{t}{N-t\gamma}\left|\mathbb{S}^{N-1}\right|
\int_0^\infty\left|u'\right|^{t-2}u'
\left(u''+\frac{N-1}{r}u'\right)
r^{N-t\gamma}\mathrm{d}r
\\&\quad+
\frac{t(N-1)}{N-t\gamma}\int_{\mathbb{R}^N}
\frac{\left|\nabla u\right|^t}{|x|^{t\gamma}}
\mathrm{d}x.
\end{align*}
This together with \eqref{2.1} implies that
\begin{align}\label{3.1}
&\frac{N-t(N+\gamma-1)}{t}
\int_{\mathbb{R}^N}
\frac{\left|\nabla u\right|^t}{|x|^{t\gamma}}
\mathrm{d}x
\nonumber\\&\qquad
=-\left|\mathbb{S}^{N-1}\right|
\int_0^\infty\left|u'\right|^{t-2}u'
\left(u''+\frac{N-1}{r}u'\right)
r^{N-t\gamma}\mathrm{d}r
\nonumber\\&\qquad
=-\left|\mathbb{S}^{N-1}\right|
\int_0^\infty\left|u'\right|^{t-2}u'
r^{-(t-1)\beta+\frac{(p-1)(N-1)}{p}}
\left(u''+\frac{N-1}{r}u'\right)
r^{-\alpha+\frac{N-1}{p}}\mathrm{d}r.
\end{align}
The density argument shows that \eqref{3.1} is still valid for all radial functions $u\in \mathcal{H}_{\alpha,\beta}^{p,\frac{p(t-1)}{p-1}}
(\mathbb{R}^N)$. By applying  H\"{o}lder inequality, we get
\begin{align}\label{3.2}
&\frac{\left|N-t(N+\gamma-1)\right|}{t}
\int_{\mathbb{R}^N}
\frac{\left|\nabla u\right|^t}{|x|^{t\gamma}}
\mathrm{d}x
\nonumber\\&\quad
\le\left(\left|\mathbb{S}^{N-1}\right|
\int_0^\infty
\left|u''+\frac{N-1}{r}u'\right|^p
r^{N-p\alpha-1}\mathrm{d}r\right)^{\frac{1}{p}}
\left[\left|\mathbb{S}^{N-1}\right|
\int_0^\infty\left|u'\right|^{\frac{p(t-1)}{p-1}}
r^{N-\frac{p(t-1)}{p-1}\beta-1}\mathrm{d}r
\right]^{\frac{p-1}{p}}
\nonumber\\&\quad
=\left(\int_{\mathbb{R}^N}
\frac{|\Delta u|^p}{|x|^{p\alpha}}
\mathrm{d}x\right)^{\frac{1}{p}}
\left[\int_{\mathbb{R}^N}
\frac{\left|\nabla u\right|^{\frac{p(t-1)}{p-1}}}
{|x|^{\frac{p(t-1)}{p-1}\beta}}
\mathrm{d}x\right]^{\frac{p-1}{p}},
\end{align}
as our desired estimate \eqref{2.2}.

Next, we provide the sufficient and necessary conditions  such that the equality in \eqref{3.2} holds. The proof is split into two cases.

$\bullet$ \emph{\textbf{Case 1:} $(p-1)(\alpha+1)-(t-1)\beta-(N-1)(t-p)>0$ and
$p\alpha+(p-1)N<0$ ($p\le t$)}.

In this case, \eqref{2.1} indicates that
\begin{align*}
N-t\left(N+\gamma-1\right)
&=-\left[p\alpha+(p-1)N\right]
+\left[(p-1)(\alpha+1)-(t-1)\beta-(N-1)(t-p)\right]>0.
\end{align*}
Thus, the equality in \eqref{3.2} holds when the H\"{o}lder inequality is used to \eqref{3.1} if and only if
$$
\left|u''+\frac{N-1}{r}u'\right|^{p-2}
\left(u''+\frac{N-1}{r}u'\right)
=-\lambda\left|u'\right|^{t-2}u'
r^{(p-1)\alpha-(t-1)\beta},
$$
for some $\lambda>0$. Let us denote $u':=r^{1-N}v$, then $v$ satisfies the following equation
$$
\left|v'\right|^{p-2}v'
=-\lambda r^{(p-1)\alpha-(t-1)\beta-(N-1)(t-p)}\left|v\right|^{t-2}v.
$$
Without losing generality, assume that $v>0$ and $v'<0$ (the proof of another case $v<0$ and $v'>0$ is the same, so we omit here). Then
$$
\left(-v'\right)^{p-1}
=\lambda r^{(p-1)\alpha-(t-1)\beta-(N-1)(t-p)}v^{t-1},
$$
which means
\begin{equation}\label{3.3}
v'
=-\lambda^{\frac{1}{p-1}}
r^{\frac{(p-1)\alpha-(t-1)\beta-(N-1)(t-p)}{p-1}}
v^{\frac{t-1}{p-1}}.
\end{equation}

$-$ If $p=t$, we have
$$
v'=-\lambda^{\frac{1}{p-1}}r^{\alpha-\beta}v,
$$
which implies that
$$
v(r)=c\mathrm{exp}
\left(-\frac{\lambda^{\frac{1}{p-1}}r^{\alpha-\beta+1}}
{\alpha-\beta+1}\right)
$$
for some $c\in\mathbb{R}$ and $\lambda>0$. Hence,
$$
u'(r)=cr^{1-N}\mathrm{exp}
\left(-\frac{\lambda^{\frac{1}{p-1}}
r^{\alpha-\beta+1}}
{\alpha-\beta+1}\right)
$$
and
$$
u(x)=c\int_{|x|}^\infty
r^{1-N}\exp\left(-\frac{\lambda^{\frac{1}{p-1}} r^{\alpha-\beta+1}}
{\alpha-\beta+1}\right)
\mathrm{d}r.
$$

$-$ If $p<t$, a simple calculation combining with \eqref{3.3} derives that
$$
\left(v^{\frac{p-t}{p-1}}\right)'
=\frac{p-t}{p-1}v^{\frac{1-t}{p-1}}v'
=\frac{t-p}{p-1}\lambda^{\frac{1}{p-1}} r^{\frac{(p-1)\alpha-(t-1)\beta-(N-1)(t-p)}{p-1}},
$$
which implies
\begin{align*}
v(r)
&=
\left[c+
\frac{\frac{t-p}{p-1}\lambda^{\frac{1}{p-1}} r^{\frac{(p-1)(\alpha+1)-(t-1)\beta-(N-1)(t-p)}{p-1}}}
{\frac{(p-1)(\alpha+1)-(t-1)\beta-(N-1)(t-p)}{p-1}}
\right]^{\frac{p-1}{p-t}}
\\&=\left[c+
\frac{(t-p)\lambda^{\frac{1}{p-1}}
r^{\frac{(p-1)(\alpha+1)-(t-1)\beta-(N-1)(t-p)}{p-1}}}
{(p-1)(\alpha+1)-(t-1)\beta-(N-1)(t-p)} \right]^{\frac{p-1}{p-t}},
\end{align*}
for some $c>0$. Therefore, the equality in \eqref{3.2} holds if and only if
\begin{align*}
u(x)
=\int_{|x|}^\infty
r^{1-N}
\left[c+
\frac{(t-p)\lambda^{\frac{1}{p-1}}
r^{\frac{(p-1)(\alpha+1)-(t-1)\beta-(N-1)(t-p)}{p-1}}}
{(p-1)(\alpha+1)-(t-1)\beta-(N-1)(t-p)}
\right]^{\frac{p-1}{p-t}}
\mathrm{d}r.
\end{align*}

$\bullet$ \emph{\textbf{Case 2:} $(p-1)(\alpha+1)-(t-1)\beta-(N-1)(t-p)<0$ and
$p\alpha+(p-1)N>0$ ($p\le t$)}.

In this case, it yields from \eqref{2.1} that
\begin{align*}
N-t\left(N+\gamma-1\right)
&=-\left[p\alpha+(p-1)N\right]
+\left[(p-1)(\alpha+1)-(t-1)\beta-(N-1)(t-p)\right]<0.
\end{align*}
Thus, the equality in \eqref{3.2} happens if and only if
$$
\left|u''+\frac{N-1}{r}u'\right|^{p-2}
\left(u''+\frac{N-1}{r}u'\right)
=\lambda\left|u'\right|^{t-2}u'
r^{(p-1)\alpha-(t-1)\beta},
$$
for some $\lambda>0$. Let $u':=r^{1-N}v$, then $v$ satisfies
$$
\left|v'\right|^{p-2}v'
=\lambda r^{(p-1)\alpha-(t-1)\beta-(N-1)(t-p)}|v|^{t-2}v.
$$
Assume that $v>0$ and $v'>0$ without loss of generality, and then
$$
\left(v'\right)^{p-1}
=\lambda r^{(p-1)\alpha-(t-1)\beta-(N-1)(t-p)}v^{t-1},
$$
that is,
$$
v'
=\lambda^{\frac{1}{p-1}} r^{\frac{(p-1)\alpha-(t-1)\beta-(N-1)(t-p)}{p-1}}
v^{\frac{t-1}{p-1}}.
$$
Using arguments similar to those of Case 1, the equality in \eqref{3.2} happens if and only if
$$
u(x)=c\int_{|x|}^\infty
r^{1-N}\exp\left(\frac{\lambda^{\frac{1}{p-1}} r^{\alpha-\beta+1}}
{\alpha-\beta+1}\right)
\mathrm{d}r,
$$
if $p=t$;
\begin{align*}
u(x)
=\int_{|x|}^\infty
r^{1-N}
\left[c+
\frac{(p-t)\lambda^{\frac{1}{p-1}} r^{\frac{(p-1)(\alpha+1)-(t-1)\beta-(N-1)(t-p)}{p-1}}}
{(p-1)(\alpha+1)-(t-1)\beta-(N-1)(t-p)}
\right]^{\frac{p-1}{p-t}}
\mathrm{d}r,
\end{align*}
if $p<t$. The proof is completed.
\end{proof}

\subsection{Second order $L^p$-CKN identities, optimal constants and optimizers for the second order $L^p$-CKN inequalities: proof of Theorems \ref{thm-2.3}, \ref{thm-2.4}, \ref{thm-2.5} and \ref{thm-2.6}}\label{3.2}

\begin{proof}[\rm\textbf{Proof of Theorem \ref{thm-2.3}}]
Due to the definition of $\mathcal{K}_p$ (see \eqref{2.3}), we see that
\begin{align}\label{3.4}
&\int_{\mathbb{R}^N}\frac{1}
{|x|^{p\alpha}}\mathcal{K}_p
\left(-|x|^{\alpha-\beta-1}x\cdot\nabla u,\Delta u\right)\mathrm{d}x
\nonumber\\&\quad=
\int_{\mathbb{R}^N}\frac{|\Delta u|^p}
{|x|^{p\alpha}}\mathrm{d}x
+(p-1)\int_{\mathbb{R}^N}\frac{|\nabla u|^p}
{|x|^{p\beta}}\mathrm{d}x
+p\int_{\mathbb{R}^N}
|x|^{-\alpha-(p-1)\beta-1}
|\nabla u|^{p-2}\left(x\cdot\nabla u\right)\Delta u\mathrm{d}x.
\end{align}
Now, it remains to estimate the last term. Indeed, applying integration by parts, one has
\begin{align}\label{3.5}
&\int_{\mathbb{R}^N}
|x|^{-\alpha-(p-1)\beta-1}
|\nabla u|^{p-2}\left(x\cdot\nabla u\right)\Delta u\mathrm{d}x
\nonumber\\&\qquad=
-\int_{\mathbb{R}^N}\nabla
\left[|x|^{-\alpha-(p-1)\beta-1}
|\nabla u|^{p-2}\left(x\cdot\nabla u\right)\right]
\cdot\nabla u\mathrm{d}x
\nonumber\\&\qquad=
\left[\alpha+(p-1)\beta\right]
\int_{\mathbb{R}^N}
\frac{|\nabla u|^p}{|x|^{\alpha+(p-1)\beta+1}}\mathrm{d}x
\nonumber\\&\qquad\quad
+\left[\alpha+(p-1)\beta+1\right]
\int_{\mathbb{R}^N}
\frac{|\nabla u|^{p-2}\left[\left(x\cdot\nabla u\right)^2-|x|^2|\nabla u|^2\right]}{|x|^{\alpha+(p-1)\beta+3}}\mathrm{d}x
\nonumber\\&\qquad\quad
-(p-1)
\int_{\mathbb{R}^N}|x|^{-\alpha-(p-1)\beta-1}
|\nabla u|^{p-2}\left(x\cdot\nabla u\right)\nabla(\nabla u)\mathrm{d}x.
\end{align}
The Divergence theorem allows us to obtain
\begin{align}\label{3.6}
&\int_{\mathbb{R}^N}|x|^{-\alpha-(p-1)\beta-1}
|\nabla u|^{p-2}\left(x\cdot\nabla u\right)\nabla(\nabla u)\mathrm{d}x
\nonumber\\&\qquad=
\frac{1}{p}
\int_{\mathbb{R}^N}|x|^{-\alpha-(p-1)\beta-1}x\cdot
\nabla\left(|\nabla u|^p\right)\mathrm{d}x
\nonumber\\&\qquad=
-\frac{1}{p}
\int_{\mathbb{R}^N}
\mathrm{div}
\left[|x|^{-\alpha-(p-1)\beta-1}x\right]
|\nabla u|^p\mathrm{d}x
\nonumber\\&\qquad=
\frac{\alpha+(p-1)\beta-N+1}{p}
\int_{\mathbb{R}^N}
\frac{|\nabla u|^p}{|x|^{\alpha+(p-1)\beta+1}}\mathrm{d}x.
\end{align}
Substituting \eqref{3.6} into \eqref{3.5}, there holds
\begin{align}\label{3.7}
&\int_{\mathbb{R}^N}
|x|^{-\alpha-(p-1)\beta-1}
|\nabla u|^{p-2}\left(x\cdot\nabla u\right)\Delta u\mathrm{d}x
\nonumber\\&\quad=
-\frac{-\alpha-(p-1)\beta-(p-1)(N-1)}{p}
\int_{\mathbb{R}^N}
\frac{|\nabla u|^p}{|x|^{\alpha+(p-1)\beta+1}}\mathrm{d}x
\nonumber\\&\qquad+
\left[\alpha+(p-1)\beta+1\right]
\int_{\mathbb{R}^N}
\frac{|\nabla u|^{p-2}\left[\left(x\cdot\nabla u\right)^2-|x|^2|\nabla u|^2\right]}{|x|^{\alpha+(p-1)\beta+3}}\mathrm{d}x.
\end{align}
Then, combining \eqref{3.4} with \eqref{3.7}, we get
\begin{align*}
&\int_{\mathbb{R}^N}\frac{1}
{|x|^{p\alpha}}\mathcal{K}_p
\left(-|x|^{\alpha-\beta-1}x\cdot\nabla u,\Delta u\right)\mathrm{d}x
\nonumber\\&\quad=
\!\int_{\mathbb{R}^N}\!\frac{|\Delta u|^p}
{|x|^{p\alpha}}\mathrm{d}x
+(p-1)\!\int_{\mathbb{R}^N}\!\frac{|\nabla u|^p}
{|x|^{p\beta}}\mathrm{d}x
-\left[-\alpha\!-\!(p-1)\beta\!-\!(p-1)(N-1)\right]
\!\int_{\mathbb{R}^N}\!
\frac{|\nabla u|^p}{|x|^{\alpha+(p-1)\beta+1}}\mathrm{d}x
\nonumber\\&\qquad
+p\left[\alpha+(p-1)\beta+1\right]
\int_{\mathbb{R}^N}
\frac{|\nabla u|^{p-2}\left[\left(x\cdot\nabla u\right)^2-|x|^2|\nabla u|^2\right]}{|x|^{\alpha+(p-1)\beta+3}}\mathrm{d}x,
\end{align*}
as our desired estimate \eqref{2.4}.

Now, it suffices to verify \eqref{2.5}. By the definition of $\mathcal{K}_p$ (see \eqref{2.3}) and \eqref{3.7}, we get
\begin{align*}
&\int_{\mathbb{R}^N}\frac{1}
{|x|^{p\alpha}}\mathcal{K}_p
\left(-\left(\frac{\int_{\mathbb{R}^N}\frac{|\Delta u|^p}
{|x|^{p\alpha}}\mathrm{d}x}
{\int_{\mathbb{R}^N}\frac{|\nabla u|^p}
{|x|^{p\beta}}\mathrm{d}x}\right)^{\frac{1}{p^2}}
|x|^{\alpha-\beta-1}x\cdot\nabla u,\left(\frac{\int_{\mathbb{R}^N}\frac{|\nabla u|^p}
{|x|^{p\beta}}\mathrm{d}x}
{\int_{\mathbb{R}^N}\frac{|\Delta u|^p}
{|x|^{p\alpha}}\mathrm{d}x}\right)^{\frac{p-1}{p^2}}\Delta u\right)\mathrm{d}x
\\&\quad
=
\left(\frac{\int_{\mathbb{R}^N}\frac{|\nabla u|^p}
{|x|^{p\beta}}\mathrm{d}x}
{\int_{\mathbb{R}^N}\frac{|\Delta u|^p}
{|x|^{p\alpha}}\mathrm{d}x}\right)^{\frac{p-1}{p}}
\int_{\mathbb{R}^N}\frac{|\Delta u|^p}
{|x|^{p\alpha}}\mathrm{d}x
+(p-1)\left(\frac{\int_{\mathbb{R}^N}\frac{|\Delta u|^p}
{|x|^{p\alpha}}\mathrm{d}x}
{\int_{\mathbb{R}^N}\frac{|\nabla u|^p}
{|x|^{p\beta}}\mathrm{d}x}\right)^{\frac{1}{p}}
\int_{\mathbb{R}^N}\frac{|\nabla u|^p}
{|x|^{p\beta}}\mathrm{d}x
\\&\qquad
+p\int_{\mathbb{R}^N}
|x|^{-\alpha-(p-1)\beta-1}
|\nabla u|^{p-2}\left(x\cdot\nabla u\right)\Delta u\mathrm{d}x
\\&\quad
=p\left(\!\int_{\mathbb{R}^N}\!\frac{|\Delta u|^p}
{|x|^{p\alpha}}\mathrm{d}x\!\right)^{\frac{1}{p}}\!
\!\left(\!\int_{\mathbb{R}^N}\!\frac{|\nabla u|^p}
{|x|^{p\beta}}\mathrm{d}x\!\right)^{\frac{p-1}{p}}
\!\!\!-\!\left[-\alpha\!-\!(p-1)\beta\!-\!(p-1)(N-1)\right]
\!\int_{\mathbb{R}^N}\!
\frac{|\nabla u|^p}{|x|^{\alpha+(p-1)\beta+1}}\mathrm{d}x
\nonumber\\&\qquad
+p\left[\alpha+(p-1)\beta+1\right]
\int_{\mathbb{R}^N}
\frac{|\nabla u|^{p-2}\left[\left(x\cdot\nabla u\right)^2-|x|^2|\nabla u|^2\right]}{|x|^{\alpha+(p-1)\beta+3}}\mathrm{d}x,
\end{align*}
which is equivalent to \eqref{2.5}. Hence, the proof has been completed.
\end{proof}

\begin{proof}[\rm\textbf{Proof of Theorem \ref{thm-2.4}}]
According to \eqref{3.7} and the definition of $\mathcal{K}_p$ (see \eqref{2.3}), we see that
\begin{align*}
&\int_{\mathbb{R}^N}\frac{1}
{|x|^{p\alpha}}\mathcal{K}_p
\left(|x|^{\alpha-\beta-1}x\cdot\nabla u,\Delta u\right)\mathrm{d}x
\\&\quad
=
\int_{\mathbb{R}^N}\frac{|\Delta u|^p}
{|x|^{p\alpha}}\mathrm{d}x
+(p-1)\int_{\mathbb{R}^N}\frac{|\nabla u|^p}
{|x|^{p\beta}}\mathrm{d}x
-p\int_{\mathbb{R}^N}
|x|^{-\alpha-(p-1)\beta-1}
|\nabla u|^{p-2}\left(x\cdot\nabla u\right)\Delta u\mathrm{d}x
\\&\quad
=
\int_{\mathbb{R}^N}\frac{|\Delta u|^p}
{|x|^{p\alpha}}\mathrm{d}x
+(p-1)\int_{\mathbb{R}^N}\frac{|\nabla u|^p}
{|x|^{p\beta}}\mathrm{d}x
\!-\!\left[\alpha\!+\!(p-1)\beta\!+\!(p-1)(N-1)\right]
\!\int_{\mathbb{R}^N}\!
\frac{|\nabla u|^p}{|x|^{\alpha+(p-1)\beta+1}}\mathrm{d}x
\nonumber\\&\qquad
-p\left[\alpha+(p-1)\beta+1\right]
\int_{\mathbb{R}^N}
\frac{|\nabla u|^{p-2}\left[\left(x\cdot\nabla u\right)^2-|x|^2|\nabla u|^2\right]}{|x|^{\alpha+(p-1)\beta+3}}\mathrm{d}x,
\end{align*}
which implies our desired estimate \eqref{2.6}. Similarly, we get
\begin{align*}
&\int_{\mathbb{R}^N}\frac{1}
{|x|^{p\alpha}}\mathcal{K}_p
\left(\left(\frac{\int_{\mathbb{R}^N}\frac{|\Delta u|^p}
{|x|^{p\alpha}}\mathrm{d}x}
{\int_{\mathbb{R}^N}\frac{|\nabla u|^p}
{|x|^{p\beta}}\mathrm{d}x}\right)^{\frac{1}{p^2}}
|x|^{\alpha-\beta-1}x\cdot\nabla u,\left(\frac{\int_{\mathbb{R}^N}\frac{|\nabla u|^p}
{|x|^{p\beta}}\mathrm{d}x}
{\int_{\mathbb{R}^N}\frac{|\Delta u|^p}
{|x|^{p\alpha}}\mathrm{d}x}\right)^{\frac{p-1}{p^2}}\Delta u\right)\mathrm{d}x
\\&\quad
=
\left(\frac{\int_{\mathbb{R}^N}\frac{|\nabla u|^p}
{|x|^{p\beta}}\mathrm{d}x}
{\int_{\mathbb{R}^N}\frac{|\Delta u|^p}
{|x|^{p\alpha}}\mathrm{d}x}\right)^{\frac{p-1}{p}}
\int_{\mathbb{R}^N}\frac{|\Delta u|^p}
{|x|^{p\alpha}}\mathrm{d}x
+(p-1)\left(\frac{\int_{\mathbb{R}^N}\frac{|\Delta u|^p}
{|x|^{p\alpha}}\mathrm{d}x}
{\int_{\mathbb{R}^N}\frac{|\nabla u|^p}
{|x|^{p\beta}}\mathrm{d}x}\right)^{\frac{1}{p}}
\int_{\mathbb{R}^N}\frac{|\nabla u|^p}
{|x|^{p\beta}}\mathrm{d}x
\\&\qquad
-p\int_{\mathbb{R}^N}
|x|^{-\alpha-(p-1)\beta-1}
|\nabla u|^{p-2}\left(x\cdot\nabla u\right)\Delta u\mathrm{d}x
\\&\quad
=p\left(\!\int_{\mathbb{R}^N}\frac{|\Delta u|^p}
{|x|^{p\alpha}}\mathrm{d}x\!\right)^{\frac{1}{p}}
\left(\!\int_{\mathbb{R}^N}\frac{|\nabla u|^p}
{|x|^{p\beta}}\mathrm{d}x\!\right)^{\frac{p-1}{p}}
\!\!-\!\left[\alpha\!+\!(p-1)\beta\!+\!(p-1)(N-1)\right]
\!\int_{\mathbb{R}^N}\!
\frac{|\nabla u|^p}{|x|^{\alpha+(p-1)\beta+1}}\mathrm{d}x
\nonumber\\&\qquad
-p\left[\alpha+(p-1)\beta+1\right]
\int_{\mathbb{R}^N}
\frac{|\nabla u|^{p-2}\left[\left(x\cdot\nabla u\right)^2-|x|^2|\nabla u|^2\right]}{|x|^{\alpha+(p-1)\beta+3}}\mathrm{d}x,
\end{align*}
which is equivalent to \eqref{2.7}. This concludes the proof.
\end{proof}

\begin{lem}
[{\rm{\cite[Lemma 1.1]{Duy22}}}]
\label{lem-3.1}
Let $N\ge1$ and $p>1$. Then
\begin{enumerate}[itemsep=0pt, topsep=2pt, parsep=1pt]

\item[(1)] $\mathcal{K}_p
    \left(\overrightarrow{X},\overrightarrow{Y}\right)
    \ge0$ for all $\overrightarrow{X},\overrightarrow{Y}\in\mathbb{R}^N$.

\item[(2)] $\mathcal{K}_p
    \left(\overrightarrow{X},\overrightarrow{Y}\right)
    =0$ if and only if $\overrightarrow{X}=\overrightarrow{Y}$.

\end{enumerate}
\end{lem}

\begin{proof}[\rm\textbf{Proof of Theorem \ref{thm-2.5}}]
We divide the proof into three cases.

\emph{\textbf{Case (I):} $\alpha+(p-1)\beta+1>0$,  $\alpha-\beta+1>0$ and $p\alpha+(p-1)N<0$.}

Since $(x\cdot\nabla u)^2\le|x|^2|\nabla u|^2$, the equality happens if and only if $u$ is radially symmetric, then the condition $\alpha+(p-1)\beta+1>0$, \eqref{2.5} and Lemma \ref{lem-3.1}-(1) indicate that
\begin{align}\label{3.22}
&\left(\int_{\mathbb{R}^N}\frac{|\Delta u|^p}
{|x|^{p\alpha}}\mathrm{d}x\right)^{\frac{1}{p}}
\left(\int_{\mathbb{R}^N}\frac{|\nabla u|^p}
{|x|^{p\beta}}\mathrm{d}x\right)^{\frac{p-1}{p}}
\!-\!\frac{\!-\alpha\!-\!(p-1)\beta\!-\!(p-1)(N-1)}{p}
\int_{\mathbb{R}^N}
\frac{|\nabla u|^p}{|x|^{\alpha+(p-1)\beta+1}}\mathrm{d}x
\nonumber\\&\quad
\ge
\left(\!\int_{\mathbb{R}^N}\frac{|\Delta u|^p}
{|x|^{p\alpha}}\mathrm{d}x\!\right)^{\frac{1}{p}}
\!\left(\!\int_{\mathbb{R}^N}\frac{|\nabla u|^p}
{|x|^{p\beta}}\mathrm{d}x\!\right)^{\frac{p-1}{p}}
\!-\!\frac{\!-\alpha\!-\!(p-1)\beta\!-\!(p-1)(N-1)}{p}
\int_{\mathbb{R}^N}
\frac{|\nabla u|^p}{|x|^{\alpha+(p-1)\beta+1}}\mathrm{d}x
\nonumber\\&\qquad
+\left[\alpha+(p-1)\beta+1\right]
\int_{\mathbb{R}^N}
\frac{|\nabla u|^{p-2}\left[\left(x\cdot\nabla u\right)^2-|x|^2|\nabla u|^2\right]}{|x|^{\alpha+(p-1)\beta+3}}\mathrm{d}x
\nonumber\\&\quad
=\frac{1}{p}\int_{\mathbb{R}^N}\frac{1}
{|x|^{p\alpha}}\mathcal{K}_p
\left(-\left(\frac{\int_{\mathbb{R}^N}\frac{|\Delta u|^p}
{|x|^{p\alpha}}\mathrm{d}x}
{\int_{\mathbb{R}^N}\frac{|\nabla u|^p}
{|x|^{p\beta}}\mathrm{d}x}\right)^{\frac{1}{p^2}}
|x|^{\alpha-\beta-1}x\cdot\nabla u,\left(\frac{\int_{\mathbb{R}^N}\frac{|\nabla u|^p}
{|x|^{p\beta}}\mathrm{d}x}
{\int_{\mathbb{R}^N}\frac{|\Delta u|^p}
{|x|^{p\alpha}}\mathrm{d}x}\right)^{\frac{p-1}{p^2}}\Delta u\right)\mathrm{d}x
\nonumber\\&\quad
\ge0.
\end{align}
That is, for all $u\in\mathcal{H}_{\alpha,\beta}^{p,p}
(\mathbb{R}^N)$,
\begin{equation}\label{3.21}
\left(\!\int_{\mathbb{R}^N}\!\frac{|\Delta u|^p}
{|x|^{p\alpha}}\mathrm{d}x\!\!\right)^{\frac{1}{p}}
\left(\!\int_{\mathbb{R}^N}\!\frac{|\nabla u|^p}
{|x|^{p\beta}}\mathrm{d}x\!\!\right)^{\frac{p-1}{p}}
\!\!\ge\frac{\!-\alpha\!-\!(p-1)\beta\!-\!(p-1)(N\!-1)}{p}
\!\int_{\mathbb{R}^N}\!
\frac{|\nabla u|^p}{|x|^{\alpha+(p-1)\beta+1}}\mathrm{d}x,
\end{equation}
as our desired estimate \eqref{2.8}. Observe that the first equality of \eqref{3.22} happens if and only if $u$ is radially symmetric, and the last equality of \eqref{3.22} happens if and only if
$$
\mathcal{K}_p
\left(-\left(\frac{\int_{\mathbb{R}^N}\frac{|\Delta u|^p}
{|x|^{p\alpha}}\mathrm{d}x}{\int_{\mathbb{R}^N}
\frac{|\nabla u|^p}
{|x|^{p\beta}}\mathrm{d}x}\right)^{\frac{1}{p^2}}
|x|^{\alpha-\beta-1}x\cdot\nabla u,\left(\frac{\int_{\mathbb{R}^N}\frac{|\nabla u|^p}
{|x|^{p\beta}}\mathrm{d}x}{\int_{\mathbb{R}^N}
\frac{|\Delta u|^p}
{|x|^{p\alpha}}\mathrm{d}x}\right)^{\frac{p-1}{p^2}}
\Delta u\right)=0.
$$
Hence, the equality in \eqref{3.21} happens if and only if $u$ is radially symmetric and satisfying
\begin{equation}\label{3.23}
\mathcal{K}_p
\left(-\left(\frac{\int_{\mathbb{R}^N}\frac{|\Delta u|^p}
{|x|^{p\alpha}}\mathrm{d}x}{\int_{\mathbb{R}^N}
\frac{|\nabla u|^p}
{|x|^{p\beta}}\mathrm{d}x}\right)^{\frac{1}{p^2}}
|x|^{\alpha-\beta-1}x\cdot\nabla u,\left(\frac{\int_{\mathbb{R}^N}\frac{|\nabla u|^p}
{|x|^{p\beta}}\mathrm{d}x}{\int_{\mathbb{R}^N}
\frac{|\Delta u|^p}
{|x|^{p\alpha}}\mathrm{d}x}\right)^{\frac{p-1}{p^2}}\Delta u\right)=0.
\end{equation}
With the help of Lemma \ref{lem-3.1}-(2), to solve \eqref{3.23} is equivalent to solve the following equation
$$
-\left(\frac{\int_{\mathbb{R}^N}\frac{|\Delta u|^p}
{|x|^{p\alpha}}\mathrm{d}x}
{\int_{\mathbb{R}^N}\frac{|\nabla u|^p}
{|x|^{p\beta}}\mathrm{d}x}\right)^{\frac{1}{p^2}}
|x|^{\alpha-\beta-1}x\cdot\nabla u=\left(\frac{\int_{\mathbb{R}^N}\frac{|\nabla u|^p}
{|x|^{p\beta}}\mathrm{d}x}
{\int_{\mathbb{R}^N}\frac{|\Delta u|^p}
{|x|^{p\alpha}}\mathrm{d}x}\right)^{\frac{p-1}{p^2}}\Delta u,
$$
that is,
$$
u''+\frac{N-1}{r}u'=-\lambda r^{\alpha-\beta}u',
$$
where $\lambda=\left(\frac{\int_{\mathbb{R}^N}\frac{|\Delta u|^p}
{|x|^{p\alpha}}\mathrm{d}x}
{\int_{\mathbb{R}^N}\frac{|\nabla u|^p}
{|x|^{p\beta}}\mathrm{d}x}\right)^{\frac{1}{p}}>0$. Let us denote $u':=r^{1-N}w$, then the above equation turns into
$$
w'=-\lambda r^{\alpha-\beta}w,
$$
this gives that
$$
w(r)=c\mathrm{exp}\left(-\frac{\lambda r^{\alpha-\beta+1}}
{\alpha-\beta+1}\right).
$$
Thus, the equality in \eqref{3.21} is only attained by
$$
u(x)=c\int_{|x|}^\infty
r^{1-N}\exp\left(-\frac{\lambda r^{\alpha-\beta+1}}
{\alpha-\beta+1}\right)
\mathrm{d}r,
$$
for $c\in\mathbb{R}$ and $\lambda>0$.

\emph{\textbf{Case (II):}
$\alpha+(p-1)\beta+1<0$, $\alpha-\beta+1<0$ and $p\alpha+(p-1)N>0$}.

Similar arguments as those of Case (I), we can obtain that
$$
\left(\int_{\mathbb{R}^N}\frac{|\Delta u|^p}
{|x|^{p\alpha}}\mathrm{d}x\right)^{\frac{1}{p}}
\left(\int_{\mathbb{R}^N}\frac{|\nabla u|^p}
{|x|^{p\beta}}\mathrm{d}x\right)^{\frac{p-1}{p}}
\!\ge\frac{\alpha+(p-1)\beta+(p-1)(N-1)}{p}
\int_{\mathbb{R}^N}
\frac{|\nabla u|^p}{|x|^{\alpha+(p-1)\beta+1}}\mathrm{d}x,
$$
for all $u\in\mathcal{H}_{\alpha,\beta}^{p,p}
(\mathbb{R}^N)$, and the equality is only valid if
$$
u(x)=c\int_{|x|}^\infty
r^{1-N}\exp\left(\frac{\lambda r^{\alpha-\beta+1}}
{\alpha-\beta+1}\right)
\mathrm{d}r,
$$
for $c\in\mathbb{R}$ and $\lambda>0$.

\emph{\textbf{Case (III):} $\alpha+(p-1)\beta+1=0$.}

In this case, it follows from Lemma \ref{lem-3.1}-(1), Theorems \ref{thm-2.3} and \ref{thm-2.4} that \eqref{2.14} holds.

Now, we verify the attainability of the constant $\frac{-\alpha-(p-1)\beta-(p-1)(N-1)}{p}$ when $\alpha-\beta+1>0$ and $p\alpha+(p-1)N<0$. For $\alpha-\beta+1>0$ and $p\alpha+(p-1)N<0$, let
$$
U(x)=\int_{|x|}^\infty
r^{1-N}\exp\left(-\frac{r^{\alpha-\beta+1}}
{\alpha-\beta+1}\right)
\mathrm{d}r.
$$
Indeed, direct calculations show that
\begin{align}\label{3.8}
\int_{\mathbb{R}^N}\frac{|\Delta U|^p}
{|x|^{p\alpha}}\mathrm{d}x
&=
\left|\mathbb{S}^{N-1}\right|
\int_0^\infty
\left|U''(r)+\frac{N-1}{r}U'(r)\right|^pr^{N-p\alpha-1}
\mathrm{d}r
\nonumber\\&=
\left|\mathbb{S}^{N-1}\right|
\int_0^\infty
\left|r^{1-N+\alpha-\beta}
\exp\left(-
\frac{r^{\alpha-\beta+1}}
{\alpha-\beta+1}\right)\right|^p
r^{N-p\alpha-1}
\mathrm{d}r
\nonumber\\&=
\left|\mathbb{S}^{N-1}\right|
\int_0^\infty
r^{N-1+p-pN-p\beta}
\exp\left(-
\frac{pr^{\alpha-\beta+1}}
{\alpha-\beta+1}\right)
\mathrm{d}r
\nonumber\\&=
\frac{\left|\mathbb{S}^{N-1}\right|}
{\alpha-\beta+1}
\int_0^\infty
r^{\frac{N+p-pN-p\beta}
{\alpha-\beta+1}-1}
\exp\left(-
\frac{pr}{\alpha-\beta+1}\right)
\mathrm{d}r
\nonumber\\&=
\frac{\left|\mathbb{S}^{N-1}\right|}{p}
\left(\frac{\alpha-\beta+1}
{p}\right)^{\frac{N+p-pN-p\beta}
{\alpha-\beta+1}-1}
\int_0^\infty
r^{\frac{N+p-pN-p\beta}
{\alpha-\beta+1}-1}
e^{-r}
\mathrm{d}r
\nonumber\\&=
\frac{\left|\mathbb{S}^{N-1}\right|}{p}
\left(\frac{\alpha-\beta+1}
{p}\right)^{\frac{N+p-pN-p\beta}
{\alpha-\beta+1}-1}
\Gamma\left(\frac{N+p-pN-p\beta}
{\alpha-\beta+1}\right)
\end{align}
and
\begin{align}\label{3.9}
\int_{\mathbb{R}^N}\frac{\left|\nabla U\right|^p}
{|x|^{p\beta}}\mathrm{d}x
&=\left|\mathbb{S}^{N-1}\right|\int_0^\infty
\left|U'(r)\right|^pr^{N-p\beta-1}
\mathrm{d}r
\nonumber\\&=
\left|\mathbb{S}^{N-1}\right|
\int_0^\infty
\left|r^{1-N}
\exp\left(-
\frac{r^{\alpha-\beta+1}}
{\alpha-\beta+1}\right)\right|^p
r^{N-p\beta-1}
\mathrm{d}r
\nonumber\\&=
\left|\mathbb{S}^{N-1}\right|
\int_0^\infty r^{N-p\beta-1+p-pN}\exp\left(-
\frac{pr^{\alpha-\beta+1}}
{\alpha-\beta+1}\right)
\mathrm{d}r
\nonumber\\&=
\frac{\left|\mathbb{S}^{N-1}\right|}
{\alpha-\beta+1}
\int_0^\infty r^{\frac{N-p\beta+p-pN}{\alpha
-\beta+1}-1}
\exp\left(-
\frac{pr}
{\alpha-\beta+1}\right)
\mathrm{d}r
\nonumber\\&=
\frac{\left|\mathbb{S}^{N-1}\right|}{p}
\left(\frac{\alpha-\beta+1}{p}
\right)^{\frac{N-p\beta+p-pN}{\alpha
-\beta+1}-1}
\Gamma\left(\frac{N-p\beta+p-pN}{\alpha
-\beta+1}\right).
\end{align}
Meanwhile,
\begin{align}\label{3.10}
\int_{\mathbb{R}^N}
\frac{\left|\nabla U\right|^p}{|x|^{\alpha+(p-1)\beta+1}}\mathrm{d}x
&=\left|\mathbb{S}^{N-1}\right|\int_0^\infty
\left|U'(r)\right|^pr^{N-\alpha-(p-1)\beta-2}
\mathrm{d}r
\nonumber\\&=
\left|\mathbb{S}^{N-1}\right|
\int_0^\infty
\left|r^{1-N}
\exp\left(-
\frac{r^{\alpha-\beta+1}}
{\alpha-\beta+1}\right)\right|^p
r^{N-\alpha-(p-1)\beta-2}
\mathrm{d}r
\nonumber\\&=
\left|\mathbb{S}^{N-1}\right|
\int_0^\infty r^{N-\alpha-(p-1)\beta-2+p-pN}\exp\left(-
\frac{pr^{\alpha-\beta+1}}
{\alpha-\beta+1}\right)
\mathrm{d}r
\nonumber\\&=
\frac{\left|\mathbb{S}^{N-1}\right|}
{\alpha-\beta+1}
\int_0^\infty r^{\frac{N-p\beta+p-pN}
{\alpha-\beta+1}-2}
\exp\left(-
\frac{pr}
{\alpha-\beta+1}\right)
\mathrm{d}r
\nonumber\\&=
\frac{\left|\mathbb{S}^{N-1}\right|}{p}
\left(\frac{\alpha-\beta+1}
{p}\right)^{\frac{N-p\beta+p-pN}
{\alpha-\beta+1}-2}
\int_0^\infty r^{\frac{N-p\beta+p-pN}
{\alpha-\beta+1}-2}
e^{-r}
\mathrm{d}r
\nonumber\\&=
\frac{\left|\mathbb{S}^{N-1}\right|}{p}
\left(\frac{\alpha-\beta+1}
{p}\right)^{\frac{N-p\beta+p-pN}
{\alpha-\beta+1}-2}
\Gamma\left(\frac{N-p\beta+p-pN}
{\alpha-\beta+1}-1\right).
\end{align}
Hence, from \eqref{3.8}$-$\eqref{3.10}, we see that
\begin{align*}
\frac{\left(\int_{\mathbb{R}^N}\frac{\left|\Delta U\right|^p}{|x|^{p\alpha}}\mathrm{d}x\right)^{\frac{1}{p}}
\left(\int_{\mathbb{R}^N}\frac{\left|\nabla U\right|^p}
{|x|^{p\beta}}\mathrm{d}x\right)^{\frac{p-1}{p}}}
{\int_{\mathbb{R}^N}
\frac{\left|\nabla U\right|^p}{|x|^{\alpha+(p-1)\beta+1}}\mathrm{d}x}
&=
\frac{\alpha-\beta+1
}{p}
\frac{
\Gamma\left(\frac{N-p\beta+p-pN}{\alpha+1
-\beta}\right)}
{\Gamma\left(\frac{N-p\beta+p-pN}
{\alpha-\beta+1}-1\right)}
\\&
=\frac{N-p\beta+p-pN-\alpha+\beta-1}{p}
\\&
=\frac{-\alpha-(p-1)\beta-(p-1)(N-1)}{p},
\end{align*}
as our desired result, and thus $\frac{-\alpha-(p-1)\beta-(p-1)(N-1)}{p}$ is the optimal constant of \eqref{2.14}.

Similarly, we prove the constant $\frac{\alpha+(p-1)\beta+(p-1)(N-1)}{p}$ when $\alpha-\beta+1<0$ and $p\alpha+(p-1)N>0$ can be attained by
$$
U(x)=\int_{|x|}^\infty
r^{1-N}\exp\left(\frac{r^{\alpha-\beta+1}}
{\alpha-\beta+1}\right)
\mathrm{d}r,
$$
and thus the constant $\frac{\alpha+(p-1)\beta+(p-1)(N-1)}{p}$ is optimal. The proof is completed.
\end{proof}

\begin{proof}[\rm\textbf{Proof of Theorem \ref{thm-2.6}}]
According to Theorems \ref{thm-2.3}, \ref{thm-2.4} and Lemma \ref{lem-3.1}, similar arguments to those of Theorem \ref{thm-2.5} can be used to derive Theorem \ref{thm-2.6}.
\end{proof}

\subsection{Optimal constants and optimizers for the weighted second order HUP: proof of Theorem \ref{thm-2.8}}\label{3.3}

\noindent We provide some preparations before proving Theorem \ref{thm-2.8}. Let
$$
y_0:=\frac{4^{\frac{2}{3}}
\left(-\sqrt[3]{3\sqrt{57}+23}
-\sqrt[3]{23-3\sqrt{57}}
+\sqrt[3]{4}\right)}{12}
\approx -0.565197717
$$
be the unique solution of $2y^3-2y^2+1=0$. Meanwhile, it is easy to check that $-8x^3-8x^2+1=0$ has three solutions,
$$
\begin{cases}
x_1=\frac{-\sqrt{5}-1}{4}\approx -0.809016994,\\
x_2=-\frac{1}{2}=-0.5,\\
x_3=\frac{\sqrt{5}-1}{4}\approx 0.309016994,
\end{cases}
$$
respectively.

\begin{lem}\label{lem-3.2}
Let $N\ge1$, under one of the following two assumptions:
\begin{enumerate}[itemsep=0pt, topsep=0pt, parsep=0pt]

\item[(1)] $\alpha+1>0$, $N+2\alpha>0$ and
$0<\frac{\alpha+1}{N+2\alpha}<\frac{\sqrt{5}-1}{4}$;

\item[(2)] $\alpha+1>0$, $N+2\alpha<0$, $\frac{\alpha+1}{N+2\alpha}
<\frac{-\sqrt{5}-1}{4}$ and $\frac{\alpha+1}{N+2\alpha+2}<\frac{\sqrt{2}-1}{2}$,

\end{enumerate}
there holds
$$
\inf\limits_{k\in\mathbb{N}\cup\{0\}}
\left[1-\frac{4k(2\alpha+2)}{(N+2k+2\alpha)^2}\right]
\frac{(N+2k+4\alpha+2)^2}{4}
=\frac{(N+4\alpha+2)^2}{4}.
$$
\end{lem}

\begin{proof}[\rm\textbf{Proof}]
Let
$$
\mathcal{F}(N,\alpha,k)
:=\left[1-\frac{4k(2\alpha+2)}{(N+2k+2\alpha)^2}\right]
\frac{(N+2k+4\alpha+2)^2}{4},
$$
and then
\begin{align*}
\mathcal{F}(N,\alpha,k)
&
=\left[\frac{1}{4}
-\frac{k(2\alpha+2)}{(N+2k+2\alpha)^2}\right]
\left[(N+2k+2\alpha)^2+(2\alpha+2)^2
+2(2\alpha+2)(N+2k+2\alpha)\right]
\\&
=\frac{(N+2k+2\alpha)^2+(2\alpha+2)^2
+2(2\alpha+2)(N+2k+2\alpha)}{4}
\\&\quad
-\frac{k(2\alpha+2)}{(N+2k+2\alpha)^2}
\left[(N+2k+2\alpha)^2+(2\alpha+2)^2
+2(2\alpha+2)(N+2k+2\alpha)\right]
\\&
=\frac{(N+2k+2\alpha)^2}{4}
+(\alpha+1)^2
+(\alpha+1)(N+2\alpha)
\!-\!\frac{k(2\alpha+2)^3}{(N+2k+2\alpha)^2}
\!-\!\frac{2k(2\alpha+2)^2}{N+2k+2\alpha}.
\end{align*}
Now, consider the function $\mathcal{F}$, for $x\in[0,\infty)$,
\begin{align*}
\mathcal{F}(N,\alpha,x)
:=\frac{(N+2x+2\alpha)^2}{4}
+(\alpha+1)^2
+(\alpha+1)(N+2\alpha)
\!-\!\frac{x(2\alpha+2)^3}{(N+2x+2\alpha)^2}
\!-\!\frac{2x(2\alpha+2)^2}{N+2x+2\alpha}.
\end{align*}
By a simple calculation, we obtain
\begin{align*}
&\frac{\partial\mathcal{F}(N,\alpha,x)}{\partial x}
\\&\quad=(N+2x+2\alpha)
-\frac{(2\alpha+2)^3}{(N+2x+2\alpha)^2}
+\frac{4x(2\alpha+2)^3}{(N+2x+2\alpha)^3}
-\frac{2(2\alpha+2)^2}{N+2x+2\alpha}
+\frac{4x(2\alpha+2)^2}{(N+2x+2\alpha)^2}
\\&\quad=\frac{y^4-(2\alpha+2)^3y
+(2y-2N-4\alpha)(2\alpha+2)^3
-2(2\alpha+2)^2y^2
+(2y-2N-4\alpha)(2\alpha+2)^2y}{y^3}
\\&\quad=\frac{y^4+(2\alpha+2)^3y
-(2N+4\alpha)(2\alpha+2)^2y
-(2N+4\alpha)(2\alpha+2)^3}{y^3}
:=\frac{\mathcal{G}(N,\alpha,y)}{y^3},
\end{align*}
where $y:=N+2x+2\alpha\ge N+2\alpha$. Then, for all $y\ge N+2\alpha\neq0$,
\begin{align}\label{3.11}
\frac{\partial\mathcal{G}(N,\alpha,y)}{\partial y}
&=4y^3+(2\alpha+2)^3
-(2N+4\alpha)(2\alpha+2)^2
\nonumber\\&\ge4(N+2\alpha)^3
+(2\alpha+2)^3
-(2N+4\alpha)(2\alpha+2)^2
\nonumber\\&=4(N+2\alpha)^3
\left[1
+2\left(\frac{\alpha+1}{N+2\alpha}\right)^3
-2\left(\frac{\alpha+1}{N+2\alpha}\right)^2\right].
\end{align}

$\bullet$ \emph{\textbf{Case 1}: $\alpha+1>0$, $N+2\alpha>0$ and
$0<\frac{\alpha+1}{N+2\alpha}<\frac{\sqrt{5}-1}{4}$.}

In this case, it follows from \eqref{3.11} that $\frac{\partial \mathcal{G}(N,\alpha,y)}{\partial y}>0$. Thus, the function $\mathcal{G}(N,\alpha,y)$ is nondecreasing with respect to $y$. Then,
\begin{align}\label{3.12}
\mathcal{G}(N,\alpha,y)
&\ge\mathcal{G}(N,\alpha,N+2\alpha)
\nonumber\\&=(N+2\alpha)^4-8(\alpha+1)^3(N+2\alpha)
-8(\alpha+1)^2(N+2\alpha)^2
\nonumber\\&=(N+2\alpha)^4
\left[1-8\left(\frac{\alpha+1}{N+2\alpha}\right)^3
-8\left(\frac{\alpha+1}{N+2\alpha}\right)^2\right]
>0,
\end{align}
for all $y>N+2\alpha>0$ and
$0<\frac{\alpha+1}{N+2\alpha}<\frac{\sqrt{5}-1}{4}$. Thus, $\mathcal{F}$ is a nondecreasing function with respect to $x$, then $\mathcal{F}(N,\alpha,x)\ge \mathcal{F}(N,\alpha,0)$ for all $x\in[0,\infty)$, this implies that
$$
\inf\limits_{k\in\mathbb{N}\cup\{0\}}
\mathcal{F}(N,\alpha,k)
=\mathcal{F}(N,\alpha,0)
=\frac{(N+4\alpha+2)^2}{4}.
$$

$\bullet$ \emph{\textbf{Case 2}: $\alpha+1>0$, $N+2\alpha<0$, $\frac{\alpha+1}{N+2\alpha}
<\frac{-\sqrt{5}-1}{4}$ and $0<\frac{\alpha+1}{N+2\alpha+2}<\frac{\sqrt{2}-1}{2}$.}

In this case, we first consider $k\in\mathbb{N}$ (not include $k=0$), then $x\in[1,\infty)$ and $y=N+2x+2\alpha\ge N+2\alpha+2>0$, we also know that $y>N+2\alpha$. Then, it yields from \eqref{3.11} that $\frac{\partial \mathcal{G}(N,\alpha,y)}{\partial y}>0$. Therefore, $\mathcal{G}(N,\alpha,y)$ is a nondecreasing function with respect to $y$. From \eqref{3.12}, we obtain that
$\mathcal{G}(N,\alpha,y)>0$ for all $y>N+2\alpha$, $N+2\alpha<0$ and $\frac{\alpha+1}{N+2\alpha}
<\frac{-\sqrt{5}-1}{4}$. Then $\mathcal{F}$ is a nondecreasing function with respect to $x$, which means that $\mathcal{F}(N,\alpha,x)
\ge\mathcal{F}(N,\alpha,1)$ for all $x\in[1,\infty)$. Thus,
$$
\min_{k\in\mathbb{N}}\mathcal{F}(N,\alpha,k)
=\mathcal{F}(N,\alpha,1)
=\left[1-\frac{4(2\alpha+2)}{(N+2\alpha+2)^2}\right]
\frac{(N+4\alpha+4)^2}{4}.
$$
Hence, with the help of $0<\frac{\alpha+1}{N+2\alpha+2}<\frac{\sqrt{2}-1}{2}$,
\begin{equation*}
\inf\limits_{k\in\mathbb{N}\cup\{0\}}
\mathcal{F}(N,\alpha,k)
=\min\left\{\mathcal{F}(N,\alpha,1),
\mathcal{F}(N,\alpha,0)\right\}
=\frac{(N+4\alpha+2)^2}{4}.
\end{equation*}
The proof is completed.
\end{proof}

\begin{rem}
Under $N\ge1$, we analyze the assumptions of Lemma \ref{lem-3.2} as follows:
\begin{align*}
\begin{cases}
\alpha+1>0,\\
N+2\alpha>0,\\
0<\frac{\alpha+1}{N+2\alpha}<\frac{\sqrt{5}-1}{4},
\end{cases}
\Longleftrightarrow
\begin{cases}
\alpha+1>0,\\
N-(\sqrt{5}-1)\alpha-(\sqrt{5}+1)>0,
\end{cases}
\end{align*}
and
\begin{equation*}
\begin{cases}
\alpha+1>0,\\
N+2\alpha<0,\\
\frac{\alpha+1}{N+2\alpha}
<\frac{-\sqrt{5}-1}{4},\\
0<\frac{\alpha+1}{N+2\alpha+2}<\frac{\sqrt{2}-1}{2}
\end{cases}
\Longleftrightarrow
\begin{cases}
N+2\alpha<0,\\
N+(\sqrt{5}+1)\alpha+(\sqrt{5}-1)>0,\\
N-2\sqrt{2}\alpha-2\sqrt{2}>0.
\end{cases}
\end{equation*}
\end{rem}

Based on the above remark, Lemma \ref{lem-3.2} can be rewritten into the following form.

\begin{cor}\label{lem-3.4}
Under one of the following two assumptions:
\begin{enumerate}[itemsep=0pt, topsep=0pt, parsep=0pt]

\item[(1)] $\alpha+1>0$ and
$N-(\sqrt{5}-1)\alpha-(\sqrt{5}+1)>0$;

\item[(2)] $N+2\alpha<0$, $N+(\sqrt{5}+1)\alpha+(\sqrt{5}-1)>0$ and $N-2\sqrt{2}\alpha-2\sqrt{2}>0$,

\end{enumerate}
then
$$
\inf\limits_{k\in\mathbb{N}\cup\{0\}}
\left[1-\frac{4k(2\alpha+2)}{(N+2k+2\alpha)^2}\right]
\frac{(N+2k+4\alpha+2)^2}{4}
=\frac{(N+4\alpha+2)^2}{4}.
$$
\end{cor}

\begin{lem}\label{lem-3.5}
If $\alpha+1<0$ and $N+4\alpha+2>0$, there holds
\begin{align*}
\inf\limits_{k\in\mathbb{N}\cup\{0\}}
\frac{(N+2k+4\alpha+2)^2}{4}
=\frac{(N+4\alpha+2)^2}{4}.
\end{align*}
\end{lem}

\begin{proof}[\rm\textbf{Proof}]
Since $N+4\alpha+2>0$ and $k\in\mathbb{N}\cup\{0\}$, it is easy to prove this lemma.
\end{proof}

Under the assumptions of Theorem \ref{thm-2.8}, we see that $N+2\alpha\neq0$ and $N+2k+2\alpha\neq0$.

\begin{lem}\label{lem-3.6}
For $N\ge1$, $k\in\mathbb{N}\cup\{0\}$, $N+2\alpha\neq0$ and $N+2k+2\alpha\neq0$, then
$$
1+\min\left\{0,\frac{-4k(2\alpha+2)}
{(N+2k+2\alpha)^2}\right\}>0.
$$
\end{lem}

\begin{proof}[\rm\textbf{Proof}]
If $\alpha+1<0$,
$$
1+\min\left\{0,\frac{-4k(2\alpha+2)}
{(N+2k+2\alpha)^2}\right\}
=1>0.
$$
If $\alpha+1>0$,
\begin{align*}
1+\min\left\{0,\frac{-4k(2\alpha+2)}
{(N+2k+2\alpha)^2}\right\}
&=1-\frac{4k(2\alpha+2)}
{(N+2k+2\alpha)^2}
\\&=\frac{N^2+4k^2+4\alpha^2+4Nk+4N\alpha-8k}
{(N+2k+2\alpha)^2}
\\&=\frac{(N+2\alpha)^2+4k(k-1)+4k(N-1)}
{(N+2k+2\alpha)^2}>0.
\end{align*}
Consequently, the result of Lemma \ref{lem-3.6} is achieved. The proof is completed.
\end{proof}

Decomposing $u$ into spherical harmonics is our main tool to prove Theorem \ref{thm-2.8}, which is a very useful method, see \cite{Cazacu22,Duong21,Tertikas07,Vazquez00} for details. In view of this, decompose a function $u\in\mathcal{C}_0^\infty(\mathbb{R}^N)$ into spherical harmonics as
$$
u(x)
=\sum_{k=0}^\infty u_k(r)\phi_k(\sigma),
$$
where $\phi_k(\sigma)$ denotes the orthonormal eigenfunctions of the Laplace-Beltrami operator $-\Delta_{\mathbb{S}^{N-1}}$ with corresponding eigenvalue $c_k=k(k+N-2)$ for $k\ge0$, that is,
$$
-\Delta_{\mathbb{S}^{N-1}}\phi_k=c_k\phi_k.
$$
The functions $u_k\in\mathcal{C}_0^\infty(\mathbb{R}^N)$ satisfy $u_k(r)=O(r^k)$, $u_k'(r)=O(r^{k-1})$ as $r\to0$, see \cite[p. 418]{Tertikas07}.

\begin{proof}[\rm\textbf{Proof of Theorem \ref{thm-2.8}}]
Let $u_k(r)=r^kw_k(r)$, from
\cite[(2.9), (2.10)]{Duong21}, we get
\begin{align}
\int_{\mathbb{R}^N}
\frac{|\Delta u|^2}{|x|^{2\alpha}}\mathrm{d}x
&=\sum_{k=0}^\infty\int_0^\infty
\left(w_k''+\frac{N+2k-1}{r}w_k'\right)^2r^{N+2k-2\alpha-1}
\mathrm{d}r,
\label{3.13}
\\
\!\int_{\mathbb{R}^N}\!|x|^{2\alpha+2}|\nabla u|^2\mathrm{d}x
&=\!\sum\limits_{k=0}^\infty
\bigg[\!\int_0^\infty\!\! r^{N+2k+2\alpha+1}\left|w_k'\right|^2\!\mathrm{d}r
\!-\!(2\alpha+2)k\!\int_0^\infty\!\! r^{N+2k+2\alpha-1}\left|w_k\right|^2\!\mathrm{d}r\bigg],
\label{3.14}\\
\int_{\mathbb{R}^N}|\nabla u|^2\mathrm{d}x
&=\sum\limits_{k=0}^\infty\int_0^\infty r^{N+2k-1}\left|w_k'\right|^2\mathrm{d}r.
\label{3.15}
\end{align}
From \eqref{3.13}, \eqref{3.14} and \eqref{3.15}, we see that \eqref{2.12} changes into
\begin{align}\label{3.16}
&
\sum\limits_{k=0}^\infty
\left[\int_0^\infty
\left(w_k''+\frac{N+2k-1}{r}w_k'\right)^2r^{N+2k-2\alpha-1}
\mathrm{d}r\right]
\nonumber\\&\quad\times
\sum\limits_{k=0}^\infty
\left[\int_0^\infty r^{N+2k+2\alpha+1}\left|w_k'\right|^2\mathrm{d}r
-(2\alpha+2)k\int_0^\infty r^{N+2k+2\alpha-1}\left|w_k\right|^2\mathrm{d}r\right]
\nonumber\\&\qquad
\ge\mathcal{S}(N,\alpha)
\left(
\sum\limits_{k=0}^\infty\int_0^\infty r^{N+2k-1}\left|w_k'\right|^2\mathrm{d}r
\right)^2.
\end{align}
Due to Cauchy-Schwarz inequality, to show \eqref{3.16}, it is enough to verify that, for all $k\ge0$,
\begin{align}\label{3.17}
&
\left[\int_0^\infty
\left(w_k''+\frac{N+2k-1}{r}w_k'\right)^2r^{N+2k-2\alpha-1}
\mathrm{d}r\right]
\nonumber\\&\quad\times
\left[\int_0^\infty r^{N+2k+2\alpha+1}\left|w_k'\right|^2\mathrm{d}r
-(2\alpha+2)k\int_0^\infty r^{N+2k+2\alpha-1}\left|w_k\right|^2\mathrm{d}r\right]
\nonumber\\&\qquad
\ge\mathcal{S}(N,\alpha)
\left(\int_0^\infty r^{N+2k-1}\left|w_k'\right|^2\mathrm{d}r
\right)^2.
\end{align}
One-dimensional Hardy inequality can be used to derive that
\begin{equation*}
\int_0^\infty r^{N+2k+2\alpha+1}\left|w_k'\right|^2\mathrm{d}r
\ge\frac{(N+2k+2\alpha)^2}{4}
\int_0^\infty r^{N+2k+2\alpha-1}\left|w_k\right|^2\mathrm{d}r.
\end{equation*}
Then, to show \eqref{3.17}, it remains to prove that, for all $k\ge0$,
\begin{align}\label{3.19}
&\int_0^\infty
\left(w_k''+\frac{N+2k-1}{r}w_k'\right)^2r^{N+2k-2\alpha-1}
\mathrm{d}r
\nonumber\\&\quad\times
\left[1+\min\left\{0,\frac{-4k(2\alpha+2)}
{(N+2k+2\alpha)^2}\right\}\right]
\int_0^\infty r^{N+2k+2\alpha+1}\left|w_k'\right|^2\mathrm{d}r
\nonumber\\&\qquad
\ge\mathcal{S}(N,\alpha)\left(\int_0^\infty r^{N+2k-1}\left|w_k'\right|^2\mathrm{d}r\right)^2.
\end{align}
Applying \cite[(1.8)]{Duong21} (with $t=2$ and $\beta=-2\alpha-2$) for radial function $w_k$ on dimension $N+2k$, we get, for all $k\ge0$,
\begin{align}\label{3.20}
&\left[\int_0^\infty
\left(w_k''+\frac{N+2k-1}{r}w_k'\right)^2r^{N+2k-2\alpha-1}
\mathrm{d}r\right]
\left[\int_0^\infty r^{N+2k+2\alpha+1}\left|w_k'\right|^2\mathrm{d}r\right]
\nonumber\\&\qquad
\ge\frac{(N+2k+4\alpha+2)^2}{4}\left(\int_0^\infty r^{N+2k-1}\left|w_k'\right|^2\mathrm{d}r\right)^2,
\end{align}
where the constant $\frac{(N+2k+4\alpha+2)^2}{4}$ of \eqref{3.20} is optimal. Thus, from Lemma \ref{lem-3.6} and \eqref{3.19},
\begin{equation*}
\mathcal{S}(N,\alpha)
\ge\inf\limits_{k\in\mathbb{N}\cup\{0\}}
\left[1+\min\left\{0,\frac{-4k(2\alpha+2)}
{(N+2k+2\alpha)^2}\right\}\right]
\frac{(N+2k+4\alpha+2)^2}{4}.
\end{equation*}
Then, with the help of Lemmas \ref{lem-3.2} and \ref{lem-3.5},
\begin{align*}
&\inf\limits_{k\in\mathbb{N}\cup\{0\}}
\left[1+\min\left\{0,\frac{-4k(2\alpha+2)}
{(N+2k+2\alpha)^2}\right\}\right]
\frac{(N+2k+4\alpha+2)^2}{4}
\\&\qquad=
\begin{cases}
\inf\limits_{k\in\mathbb{N}\cup\{0\}}
\left[1-\frac{4k(2\alpha+2)}{(N+2k+2\alpha)^2}\right]
\frac{(N+2k+4\alpha+2)^2}{4},
\ \ &\ \mathrm{if}\ \alpha+1>0\\
\inf\limits_{k\in\mathbb{N}\cup\{0\}}
\frac{(N+2k+4\alpha+2)^2}{4},
\ \ &\ \mathrm{if}\ \alpha+1<0
\end{cases}
\\&\qquad
=\frac{(N+4\alpha+2)^2}{4}.
\end{align*}
Hence,
$$
\int_{\mathbb{R}^N}
\frac{|\Delta u|^2}{|x|^{2\alpha}}\mathrm{d}x
\int_{\mathbb{R}^N}
|x|^{2\alpha+2}|\nabla u|^2\mathrm{d}x
\ge\frac{\left(N+4\alpha+2\right)^2}{4}
\left(\int_{\mathbb{R}^N}|\nabla u|^2\mathrm{d}x\right)^2,
$$
equivalently,
\begin{equation*}
\left(\int_{\mathbb{R}^N}
\frac{|\Delta u|^2}{|x|^{2\alpha}}\mathrm{d}x
\right)^{\frac{1}{2}}
\left(\int_{\mathbb{R}^N}
|x|^{2\alpha+2}|\nabla u|^2\mathrm{d}x
\right)^{\frac{1}{2}}
\ge\frac{|N+4\alpha+2|}{2}
\int_{\mathbb{R}^N}|\nabla u|^2\mathrm{d}x.
\end{equation*}

Next, we verify the attainability of the constant $\frac{N+4\alpha+2}{2}$ when $\alpha+1>0$ and $N+2\alpha>0$. For $\alpha+1>0$ and $N+2\alpha>0$, let
$$
U(x)=e^{-\frac{|x|^{2(\alpha+1)}}{2(\alpha+1)}}.
$$
Indeed, direct calculations yield that
\begin{align*}
\int_{\mathbb{R}^N}
\frac{|\Delta U|^2}{|x|^{2\alpha}}\mathrm{d}x
&=\left|\mathbb{S}^{N-1}\right|
\int_0^\infty\left|U''(r)+\frac{N-1}{r}U'(r)\right|^2
r^{N-2\alpha-1}\mathrm{d}r
\\&=\left|\mathbb{S}^{N-1}\right|
\int_0^\infty
\left|r^{2\alpha+2}-(N+2\alpha)\right|^2
e^{-\frac{r^{2(\alpha+1)}}{\alpha+1}}
r^{N+2\alpha-1}\mathrm{d}r
\\&=\left|\mathbb{S}^{N-1}\right|
\bigg[\int_0^\infty
e^{-\frac{r^{2(\alpha+1)}}{\alpha+1}}
r^{N+6\alpha+3}\mathrm{d}r
-2(N+2\alpha)\int_0^\infty
e^{-\frac{r^{2(\alpha+1)}}{\alpha+1}}
r^{N+4\alpha+1}\mathrm{d}r
\\&\quad
+(N+2\alpha)^2\int_0^\infty
e^{-\frac{r^{2(\alpha+1)}}{\alpha+1}}
r^{N+2\alpha-1}\mathrm{d}r\bigg],
\end{align*}
where
\begin{align*}
\int_0^\infty
e^{-\frac{r^{2(\alpha+1)}}{\alpha+1}}
r^{N+6\alpha+3}\mathrm{d}r
&=\frac{1}{2}(\alpha+1)^{\frac{N+4\alpha+2}{2\alpha+2}}
\Gamma\left(\frac{N+6\alpha+4}{2\alpha+2}\right),\\
\int_0^\infty
e^{-\frac{r^{2(\alpha+1)}}{\alpha+1}}
r^{N+4\alpha+1}\mathrm{d}r
&=\frac{1}{2}(\alpha+1)^{\frac{N+2\alpha}{2\alpha+2}}
\Gamma\left(\frac{N+4\alpha+2}{2\alpha+2}\right),\\
\int_0^\infty
e^{-\frac{r^{2(\alpha+1)}}{\alpha+1}}
r^{N+2\alpha-1}\mathrm{d}r
&=\frac{1}{2}(\alpha+1)^{\frac{N-2}{2\alpha+2}}
\Gamma\left(\frac{N+2\alpha}{2\alpha+2}\right).
\end{align*}
Hence,
\begin{equation*}
\int_{\mathbb{R}^N}
\frac{|\Delta U|^2}{|x|^{2\alpha}}\mathrm{d}x
=\frac{1}{8}\left|\mathbb{S}^{N-1}\right|
(N+2\alpha)(N+4\alpha+2)
(\alpha+1)^{\frac{N-2}{2\alpha+2}}
\Gamma\left(\frac{N+2\alpha}{2\alpha+2}\right).
\end{equation*}
Meanwhile, there hold
\begin{align*}
\int_{\mathbb{R}^N}
|x|^{2\alpha+2}|\nabla U|^2\mathrm{d}x
&=\left|\mathbb{S}^{N-1}\right|
\int_0^\infty\left|U'(r)\right|^2
r^{N+2\alpha+1}\mathrm{d}r
\\&=\left|\mathbb{S}^{N-1}\right|
\int_0^\infty e^{-\frac{r^{2(\alpha+1)}}{\alpha+1}}
r^{N+6\alpha+3}\mathrm{d}r
\\&=\frac{1}{2}\left|\mathbb{S}^{N-1}\right|
(\alpha+1)^{\frac{N+4\alpha+2}{2\alpha+2}}
\Gamma\left(\frac{N+6\alpha+4}{2\alpha+2}\right)
\end{align*}
and
\begin{align*}
\int_{\mathbb{R}^N}
|\nabla U|^2\mathrm{d}x
&=\left|\mathbb{S}^{N-1}\right|
\int_0^\infty\left|U'(r)\right|^2
r^{N-1}\mathrm{d}r
\\&=\left|\mathbb{S}^{N-1}\right|
\int_0^\infty e^{-\frac{r^{2(\alpha+1)}}{\alpha+1}}
r^{N+4\alpha+1}\mathrm{d}r
\\&=\frac{1}{2}\left|\mathbb{S}^{N-1}\right|
(\alpha+1)^{\frac{N+2\alpha}{2\alpha+2}}
\Gamma\left(\frac{N+4\alpha+2}{2\alpha+2}\right).
\end{align*}
Therefore,
\begin{equation*}
\frac{\left(\int_{\mathbb{R}^N}
\frac{|\Delta U|^2}{|x|^{2\alpha}}\mathrm{d}x\right)^{\frac{1}{2}}
\left(\int_{\mathbb{R}^N}
|x|^{2\alpha+2}|\nabla U|^2\mathrm{d}x\right)^{\frac{1}{2}}}
{\int_{\mathbb{R}^N}|\nabla U|^2\mathrm{d}x}
=\frac{N+4\alpha+2}{2},
\end{equation*}
as our desired, and thus $\frac{N+4\alpha+2}{2}$ is the  optimal constant of \eqref{2.12}.

We can similarly verify that the constant $-\frac{N+4\alpha+2}{2}$ when $\alpha+1<0$ and $N+2\alpha<0$ can be attained by
$$
U(x)=e^{\frac{|x|^{2(\alpha+1)}}{2(\alpha+1)}},
$$
and thus $-\frac{N+4\alpha+2}{2}$ is the  optimal constant of \eqref{2.12}. The proof is completed.
\end{proof}

\begin{rem}
We analyze the assumptions of Theorem \ref{thm-2.8} as follows.

\begin{enumerate}[itemsep=0pt, topsep=0pt, parsep=0pt]

\item[(1)] We apply \cite[(1.8)]{Duong21} (with $t=2$ and $\beta=-2\alpha-2$) to obtain \eqref{3.20} above, there need
$$
N-2\alpha>0,
\
N+2\alpha+2>0.
$$
This together with the assumptions of Corollary \ref{lem-3.4} and Lemma \ref{lem-3.5}, the assumptions of Theorem \ref{thm-2.8} are obtained.

\item[(2)] In order to verify the attainability of the optimal constants $\pm\frac{N+4\alpha+2}{2}$, there need to apply $\frac{N+2\alpha}{\alpha+1}>0$ to guarantee that $\Gamma\left(\frac{N+6\alpha+4}{2\alpha+2}\right)$, $\Gamma\left(\frac{N+4\alpha+2}{2\alpha+2}\right)$ and $\Gamma\left(\frac{N+2\alpha}{2\alpha+2}\right)$ make sense.

\end{enumerate}
\end{rem}

\section*{disclosure statement}
\noindent No potential conflict of interest was reported by the author(s).

\section*{Data availability statements}
\noindent No data was used for the research described in the article.

\section*{Acknowledgements}
\noindent This paper was supported by the National Natural Science Foundation of China (No. 12371120).


\begin{thebibliography}{99}
\bibitem{Brezis18}
H. Brezis, P. Mironescu,
Gagliardo-Nirenberg inequalities and non-inequalities: the full story.
\emph{Ann. Inst. H. Poincar\'{e} Anal. Non Lin\'{e}aire} \textbf{35} (2018), no. 5, 1355--1376.


\bibitem{Caffarelli84}
L. Caffarelli, R. Kohn, L. Nirenberg,
First order interpolation inequalities with weights. \emph{Compositio Math.} \textbf{53} (1984), no. 3, 259--275.


\bibitem{Catrina09}
F. Catrina, D.G. Costa,
Sharp weighted-norm inequalities for functions with compact support in $\mathbb{R}^N\setminus\{0\}$.
\emph{J. Differential Equations} \textbf{246} (2009), no. 1, 164--182.


\bibitem{Catrina01}
F. Catrina, Z.-Q. Wang,
On the Caffarelli-Kohn-Nirenberg inequalities: sharp constants, existence (and nonexistence), and symmetry of extremal functions.
\emph{Comm. Pure Appl. Math.} \textbf{54} (2001), no. 2, 229--258.


\bibitem{Cazacu20}
C. Cazacu,
A new proof of the Hardy-Rellich inequality in any dimension.
\emph{Proc. Roy. Soc. Edinburgh Sect. A} \textbf{150} (2020), no. 6, 2894--2904.


\bibitem{Cazacu21}
C. Cazacu, J. Flynn, N. Lam,
Short proofs of refined sharp Caffarelli-Kohn-Nirenberg inequalities.
\emph{J. Differential Equations} \textbf{302} (2021), 533--549.


\bibitem{Cazacu22}
C. Cazacu, J. Flynn, N. Lam,
Sharp second order uncertainty principles.
\emph{J. Funct. Anal.} \textbf{283} (2022), no. 10, Paper No. 109659, 37 pp.


\bibitem{Cazacu23}
C. Cazacu, J. Flynn, N. Lam,
Caffarelli-Kohn-Nirenberg inequalities for curl-free vector fields and second order derivatives.
\emph{Calc. Var. Partial Differential Equations} \textbf{62} (2023), no. 4, Paper No. 118, 26 pp.


\bibitem{Cazacu24}
C. Cazacu, D. Krej\v{c}i\v{r}\'{i}k, N. Lam, A. Laptev,
Hardy inequalities for magnetic $p$-Laplacians. \emph{Nonlinearity} \textbf{37} (2024), no. 3, Paper No. 035004, 27 pp.


\bibitem{Cordero-Erausquin04}
D. Cordero-Erausquin, B. Nazaret, C. Villani,
A mass-transportation approach to sharp Sobolev and Gagliardo-Nirenberg inequalities.
\emph{Adv. Math.} \textbf{182} (2004), no. 2, 307--332.


\bibitem{Costa08}
D.G. Costa,
Some new and short proofs for a class of Caffarelli-Kohn-Nirenberg type inequalities.
\emph{J. Math. Anal. Appl.} \textbf{337} (2008), no. 1, 311--317.


\bibitem{DelPino02}
M. Del Pino, J. Dolbeault,
Best constants for Gagliardo-Nirenberg inequalities and applications to nonlinear diffusions.
\emph{J. Math. Pures Appl. (9)} \textbf{81} (2002), no. 9, 847--875.


\bibitem{DelPino03}
M. Del Pino, J. Dolbeault,
The optimal Euclidean $L^p$-Sobolev logarithmic inequality.
\emph{J. Funct. Anal.} \textbf{197} (2003), no. 1, 151--161.


\bibitem{Do23}
A.X. Do, J. Flynn, N. Lam, G. Lu,
$L^{p}$-Caffarelli-Kohn-Nirenberg inequalities and their stabilities. arXiv preprint arXiv:2310.07083.


\bibitem{Do23-2}
A.X. Do, N. Lam, G. Lu,
A new approach to weighted Hardy-Rellich inequalities: improvements, symmetrization principle and symmetry breaking. arXiv preprint arXiv:2310.06637.


\bibitem{Dong18}
M. Dong, N. Lam, G. Lu,
Sharp weighted Trudinger-Moser and Caffarelli-Kohn-Nirenberg inequalities and their extremal functions.
\emph{Nonlinear Anal.} \textbf{173} (2018), 75--98.


\bibitem{Duong21}
A.T. Duong, V.H. Nguyen,
The sharp second order Caffareli-Kohn-Nirenberg inequality and stability estimates for the sharp second order uncertainty principle. arXiv preprint arXiv:2102.01425.


\bibitem{Duy22}
N.T. Duy, N. Lam, G. Lu,
$p$-Bessel pairs, Hardy's identities and inequalities and Hardy-Sobolev inequalities with monomial weights.
\emph{J. Geom. Anal.} \textbf{32} (2022), no. 4, Paper No. 109, 36 pp.


\bibitem{Folland97}
G.B. Folland, A. Sitaram,
The uncertainty principle: a mathematical survey.
\emph{J. Fourier Anal. Appl.} \textbf{3} (1997), no. 3, 207--238.


\bibitem{Frank08}
R.L. Frank, R. Seiringer,
Non-linear ground state representations and sharp Hardy inequalities.
\emph{J. Funct. Anal.} \textbf{255} (2008), no. 12, 3407--3430.


\bibitem{Lam20}
N. Lam, G. Lu, L. Zhang,
Geometric Hardy's inequalities with general distance functions.
\emph{J. Funct. Anal.} \textbf{279} (2020), no. 8, 108673, 35 pp.


\bibitem{Lieb10}
E.H. Lieb, R. Seiringer,
\emph{The stability of matter in quantum mechanics}.
Cambridge University Press, Cambridge, 2010.


\bibitem{Lin86}
C.S. Lin,
Interpolation inequalities with weights.
\emph{Comm. Partial Differential Equations} \textbf{11} (1986), no. 14, 1515--1538.


\bibitem{Nguyen15}
V.H. Nguyen,
Sharp weighted Sobolev and Gagliardo-Nirenberg inequalities on half-spaces via mass transport and consequences.
\emph{Proc. Lond. Math. Soc. (3)} \textbf{111} (2015), no. 1, 127--148.


\bibitem{Talenti76}
G. Talenti,
Best constant in Sobolev inequality.
\emph{Ann. Mat. Pura Appl. (4)} \textbf{110} (1976), 353--372.


\bibitem{Tertikas07}
A. Tertikas, N.B. Zographopoulos,
Best constants in the Hardy-Rellich inequalities and related improvements.
\emph{Adv. Math.} \textbf{209} (2007), no. 2, 407--459.


\bibitem{Vazquez00}
J.L. Vazquez, E. Zuazua,
The Hardy inequality and the asymptotic behaviour of the heat equation with an inverse-square potential.
\emph{J. Funct. Anal.} \textbf{173} (2000), no. 1, 103--153.


\bibitem{Weyl50}
H. Weyl,
\emph{The theory of groups and quantum mechanics}. Translated from the second (revised) German edition by H. P. Robertson. Reprint of the 1931 English translation. Dover Publications, Inc., New York, 1950.


\bibitem{Xia07}
C. Xia,
The Caffarelli-Kohn-Nirenberg inequalities on complete manifolds.
\emph{Math. Res. Lett.} \textbf{14} (2007), no. 5, 875--885.


\end{thebibliography}
\end{document}